\documentclass[11pt]{amsart}
\usepackage{amsbsy,amssymb,amsmath,amsthm,amscd,amsfonts,latexsym,amstext,delarray,
amsmath,graphicx} 
\usepackage[margin=1in]{geometry}
\usepackage{color}
\input xypic

\usepackage[all]{xy}
\usepackage{setspace}

\newtheorem{thm}{Theorem}[section]
\newtheorem{prop}[thm]{Proposition}
\newtheorem{cor}[thm]{Corollary}
\newtheorem{lem}[thm]{Lemma}

\newtheorem{defn}[thm]{Definition}
\newtheorem{rem}[thm]{Remark}
\newtheorem{ex}[thm]{Example}

\numberwithin{equation}{section}

\def\F{{\mathbb F}}
\def\Q{{\mathbb Q}}
\def\Z{{\mathbb Z}}
\def\N{{\mathbb N}}
\def\R{{\mathbb R}}
\def\C{{\mathbb C}}

\def\A{{\mathbb A}}

\def\cA{{\mathcal A}}

\def\cC{{\mathcal C}}
\def\cD{{\mathcal D}}

\def\cG{{\mathcal G}}
\def\cH{{\mathcal H}}

\def\cL{{\mathcal L}}
\def\cM{{\mathcal M}}
\def\cN{{\mathcal N}}

\def\cP{{\mathcal P}}

\def\cR{{\mathcal R}}
\def\cS{{\mathcal S}}
\def\cT{{\mathcal T}}

\def\cW{{\mathcal W}}

\def\bL{{\mathbb L}}

\def\bR{{\mathbb R}}
\def\bS{{\mathbb S}}
\def\bT{{\mathbb T}}

\def\fA{{\mathfrak A}}

\def\Tr{{\rm Tr}}
\def\Spec{{\rm Spec}}
\def\m{{\mathfrak m}}

\DeclareMathOperator*{\Hom}{Hom}

\title{Entropy algebras and Birkhoff factorization}
\author{Matilde Marcolli and Nicolas Tedeschi}
\address{Mathematics Department, Caltech, 1200 E. California Blvd. Pasadena, CA 91125, USA}
\email{matilde@caltech.edu}
\email{nicot@caltech.edu}
\date{}

\begin{document}
\maketitle

\begin{abstract}
We develop notions of Rota--Baxter structures and associated Birkhoff factorizations,
in the context of min-plus semirings and their thermodynamic deformations, including
deformations arising from quantum information measures such as the von Neumann
entropy. We consider examples related to Manin's renormalization and computation 
program, to Markov random fields and to counting functions and zeta functions of
algebraic varieties.
\end{abstract}

\section{Introduction}

This paper is motivated by two different sources: Manin's ``renormalization and computation"
program, \cite{Man1}, \cite{Man2}, \cite{Man3}, and the theory of ``thermodynamic semirings"
developed in \cite{CoCo}, \cite{MaThor}. Manin proposed the use of an algebraic framework
modeled on the Connes--Kreimer theory of renormalization \cite{CoKr} to achieve a renormalization
of infinities that arise in computation (halting problem). In this formalism the Hopf algebra of Feynman
graphs is replaced by a Hopf algebra of flow charts computing recursive functions. He suggested that
natural characters of this Hopf algebra, of relevance to the computational setting, such as memory
size or computing time, would be naturally taking values in a min-plus or max-plus (tropical) algebra instead
of taking values in a commutative Rota--Baxter algebra, as in the case of renormalization in
quantum field theory. Thus, in \cite{Man1} he asked for an extension of the algebraic
renormalization method based on Rota--Baxter algebras (\cite{EFGuo}, \cite{EFGK}) to
tropical semirings. On the other hand, min-plus semirings admit deformations based on
thermodynamic information measures, such as the Shannon entropy and generalizations,
\cite{MaThor}. These are closely related to Maslov dequantization, \cite{Viro}.

\smallskip

In this paper we develop a unified approach to Rota--Baxter structures and 
Birkhoff factorizations in min-plus semirings and their thermodynamic deformations.
In Section \S \ref{thermoSec} we recall the basic definitions and properties of
thermodynamic semirings, and we describe generalizations defined as deformations
of the (tropical) trace using the von Neumann entropy and other entropy measures in
quantum information. In \S \ref{RBBirkSec} we introduce Rota--Baxter structures
on min-plus semirings and we obtain a Birkhoff factorization of min-plus characters
for Rota--Baxter structures of weight $+1$ (unlike the original renormalization
case that uses weight $-1$ Rota--Baxter operators). In \S \ref{thermoRBsec}, we
introduce Rota--Baxter structures on thermodynamic semirings and we relate them
to Rota--Baxter structures on ordinary commutative rings. We construct 
Birkhoff factorizations in thermodynamic semirings with Rota--Baxter
operators of weight $+1$. In \S \ref{vonNeumannSec} we extend 
the thermodynamic Rota--Baxter structures to the case of the von Neumann 
entropy and the trace deformation. In \S \ref{WittSec} we consider some
explicit examples of Rota--Baxter operators of weight $+1$ on commutative
rings and on thermodynamic semirings, and we show that they determine
Rota--Baxter structures of the same weight on Witt rings. We discuss some
applications to zeta functions of algebraic varieties, seen as elements of
Witt rings, as in \cite{Rama}. In \S \ref{minoneSec} we consider Rota--Baxter
operators of weight $-1$ on min-plus semirings, and we show that, under an
additional superadditivity condition, one can still obtain Birkhoff factorizations.
We also consider a variant of the construction, where the Birkhoff factorization
is obtained from a pair of Rota--Baxter operators of weight $-1$, generalizing
the pair $T$, $id-T$ of the classical renormalization case. In \S \ref{examplesSec}
we consider three explicit examples of min-plus characters, motivated, respectively,
by Manin's renormalization and computation proposal \cite{Man1}, and the complexity
theory of recursive functions \cite{Blum}; by the theory of Markov random fields
and Gibbs states on graphs, \cite{Pres}; and by the question of polynomial
countability for the graph hypersurfaces of quantum field theory, \cite{Mar}.

\bigskip
\section{Thermodynamic semirings and other thermodynamic deformations}\label{thermoSec}

After recalling the notion of thermodynamic semirings from \cite{MaThor}, we 
introduce thermodynamic deformations of the trace, which extend the deformed
addition of thermodynamic semirings from classical to quantum information.
In particular, we interpret the case based on the von Neumann entropy as a Helmholtz free energy.
We also discuss briefly functionals obtained as thermodynamic deformations
of the integral, defined through the data of a dynamical system and its metric 
and topological entropies.

\subsection{Thermodynamic semirings}

The min-plus (or tropical) semiring $\bT$ is $\bT=\R \cup \{ \infty \}$, with the
operations $\oplus$ and $\odot$ given by
$$ x \oplus y = \min \{ x, y \}, $$
with $\infty$ the identity element for $\oplus$ and with
$$ x \odot y = x + y, $$
with $0$ the identity element for $\odot$. The operations $\oplus$ and $\odot$
satisfy associativity and commutativity and distributivity of the product $\odot$ 
over the sum $\oplus$. 

\smallskip

We will occasionally consider also the analogous max-plus version $\bT_{\rm max}=\R\cup\{ -\infty \}$,
with $\oplus=\max$ and $\odot=+$. We will write $\bT_{\rm max}$, when needed, to 
distinguish it from $\bT=\bT_{\rm min}$.

\smallskip

A notion of {\em thermodynamic semiring} was developed in \cite{MaThor}, generalizing
a construction of \cite{CoCo}, 
as a deformation of the min-plus algebra, where the product $\odot$ is
unchanged, but the sum $\oplus$ is deformed to a new operation $\oplus_{\beta,S}$,
according to a binary entropy functional $S$ and a deformation parameter $\beta\geq 0$,
which we interpret thermodynamically as an inverse temperature  (up to
the Boltzmann constant which we set equal to $1$). At zero temperature (that is,
$\beta\to \infty$) one recovers the unperturbed idempotent addition.
The case where the entropy functional $S$ is the Shannon entropy was
considered in \cite{CoCo}, in relation to geometry
over the field with one element, while other entropy functionals, such as
R\'enyi entropy or Tsallis entropy or Kullback--Leibler divergence are
considered in \cite{MaThor}, along with a general operadic formulation.

\smallskip

More precisely, for a fixed $\beta \geq 0$ and a given entropy functional $S$, 
one defines on $\R\cup \{ \infty \}$ the operation
\begin{equation}\label{oplusTS}
x \oplus_{\beta,S} y = \min_p \{ p x + (1-p) y - \frac{1}{\beta} S(p) \}.
\end{equation}
The algebraic properties (commutativity, left and right identity, associativity)
of this operation correspond to properties of the entropy functional
(symmetry $S(p)=S(1-p)$, minima $S(0)=S(1)=0$, and extensivity $S(pq)+
(1-pq) S(p(1-q)/(1-pq))= S(p) +p S(q)$). Thus, by the Khinchin axioms,
imposing that all the algebraic properties of $\bT$ are preserved in the deformation singles out the
Shannon entropy among the possible functionals $S$, while non-extensive
entropy (see \cite{GeTsa}) can be modeled by non-associative thermodynamic semirings.
We refer the reader to \cite{MaThor} for more details.

\smallskip

When $S$ is the Shannon entropy,
the idempotent property $x\oplus x=\min\{ x,x \}=x$ of the tropical addition
becomes in the deformed case $x\oplus_{\beta,S} x= x - \beta^{-1} \log2$.
This is immediately evident from $x\oplus_{\beta,S} y=\min_p \{ p x + (1-p) y - \beta^{-1} S(p) \}$,
which for $y=x$ gives $x\oplus_{\beta,S} x= x - \beta^{-1} \max_p S(p) = x - \beta^{-1}\log 2$.
Moreover, in the case of the Shannon entropy, the deformed addition can be
written equivalently as
\begin{equation}\label{logexpbetaS}
x \oplus_{\beta,S} y = -\beta^{-1} \log \left( e^{-\beta x} + e^{-\beta y} \right).
\end{equation}

\smallskip

The theory of thermodynamic semirings developed in \cite{MaThor} leads 
to a more general operadic and categorical formulation of entropy functionals
(see \S 10 of \cite{MaThor}), which is similar in spirit to the approach of \cite{BaFriLei}.

\smallskip

As in \S 10 of \cite{MaThor}, consider a collection $\cS=\{ S_n \}_{n\in \N}$
of $n$-ary entropy functionals $S_n$, satisfying the coherence
condition
$$ S_n(p_1,\ldots,p_n)=S_m(p_{i_1},\ldots,p_{i_m}), $$
whenever,  for some $m<n$, we have $p_j=0$ for all 
$j\notin \{ i_1, \ldots, i_m \}$. 

\smallskip

Shannon, R\'enyi, Tsallis entropies satisfy the coherence condition, 
and so do, more generally, entropy functionals depending on functions 
$f$ and $g$ of the form
$$ S_n(p_1,\ldots, p_n)= f(\sum_{i=1}^n g(p_i)). $$

\smallskip

A collection $\cS=\{ S_n \}_{n\in \N}$ as above determines a family of $n$-ary
operations $C_{n,\beta,\cS}$ on $\R\cup \{ \infty \}$, 
\begin{equation}\label{CnTS}
C_{n,\beta,\cS}(x_1,\ldots,x_n)= \min_{p} \{ \sum_{i=1}^n p_i x_i - \frac{1}{\beta} S_n(p_1,\ldots,p_n) \},
\end{equation}
where the minimum is taken over $p=(p_i)$, with $\sum_i p_i=1$. More generally,
given $\cS$ as above and the collection of all rooted tree $\cT$ with $n$ leaves, and
with fixed planar embeddings, we obtain $n$-ary operations $C_{n,\beta,\cS,\cT}(x_1,\ldots,x_n)$
on $\R\cup \{ \infty \}$, determined by the tree $\cT$ and the collection of
entropy functionals $S_j$ for $j=2,\ldots,n+1$. Namely, one defines 
$C_{n,\beta,\cS,\cT}(x_1,\ldots,x_n)$ as the output of the tree $\cT$
with inputs $x_1, \ldots, x_n$ at the leaves and with an operation $C_{m,\beta,\cS}$
at each vertex of valence $m+1$. As shown in Theorem 10.9 of \cite{MaThor}, these
operations can be written equivalently as
\begin{equation}\label{CnTStree}
C_{n,\beta,\cS,\cT}(x_1,\ldots,x_n)= \min_{p} \{ \sum_{i=1}^n p_i x_i - \frac{1}{\beta} S_\cT(p_1,\ldots,p_n) \},
\end{equation}
with the $S_\cT(p_1,\ldots,p_n)$ obtained from the $S_j$, for $j=2,\ldots,n$.

\smallskip

The data $(\bT,\cS)$ with $\bT=(\R \cup \{ \infty \},\oplus,\odot)$ and with $\cS=\{ S_n \}_{n\in \N}$
a coherent family of entropy functionals define an {\em information algebra}, which is an algebra
over the $A_\infty$-operad of rooted trees, see \S 10 of \cite{MaThor}.

\medskip
\subsection{Von Neumann entropy and deformed traces}

When passing from classical to quantum infomation, probabilities $P=(p_i)_{i=1}^n$
with $p_i\geq 0$ and $\sum_i p_i=1$ are replaced by density matrices $\rho$ with
$\rho^*=\rho$, $\rho\geq 0$, and $\Tr(\rho)=1$. The classical case is recovered as
the case of diagonal matrices. Correspondingly, the entropy functionals, such as
Shannon entropy, R\'enyi and Tsallis entropies, Kullback--Leibler relative entropy,
have quantum information analogs, given by the von Neumann entropy and its
generalizations. The algebraic structure of thermodynamic semirings, which encodes
the axiomatic properties of classical entropy functionals, also generalizes to
quantum information, no longer in the form of a deformed addition on a semiring,
but as a deformed trace, as we discuss below.

\smallskip

For $N\geq 1$, let $$\cM^{(N)}=\{ \rho \in M_{N\times N}(\C) \,|\, \rho^*=\rho, \, \rho\geq 0, \, \Tr(\rho)=1\}$$
be the convex set of density matrices. The von Neumann entropy
\begin{equation}\label{vNent}
\cN(\rho)= -\Tr(\rho \log\rho), \ \ \ \text{ for } \ \ \rho\in \cM^{(N)},
\end{equation}
is the natural generalization of the Shannon entropy to the quantum information setting. It reduces
to the Shannon entropy in the diagonal case.

\smallskip

As above, let $\bT=(\bR\cup \{ \infty \}, \oplus, \odot)$ be the tropical min-plus semiring. 
Let $M_{N\times N}(\bT)$ denote $N\times N$-matrices with entries in $\bR\cup \{ \infty \}$,
with the operations of idempotent matrix addition and multiplication
$$ (A \oplus B)_{ij} = \min\{ A_{ij}, B_{ij} \}, \ \ \  (A\odot B)_{ij} = \oplus_k  A_{ik}\odot B_{kj}
= \min_k  \{  A_{ik} + B_{kj} \} . $$
The trace is defined as:
\begin{equation}\label{Trmin}
\Tr^\oplus(A)=\min_i \{ A_{ii} \}. 
\end{equation}
We also denote by
\begin{equation}\label{tildeTr}
\widetilde\Tr^\oplus(A) := \min_{U\in U(N)} \min_i \{ (UAU^*)_{ii} \} \leq \Tr^\oplus(A).
\end{equation}

\smallskip

We introduce thermodynamic deformations of the trace, by setting
\begin{equation}\label{SrhoTr}
\Tr^\oplus_{\beta,S}(A):=\min_{\rho \in \cM^{(N)}} \{ \Tr(\rho A) - \beta^{-1} S(\rho) \},
\end{equation}
where $\Tr$ in the right-hand-side is the ordinary trace, with $\Tr(\rho A)=\langle A \rangle$ the
expectation value of the observable $A$ with respect to the state $\varphi(\cdot) = \Tr(\rho\, \cdot)$.

\begin{lem}\label{Trbeta0}
The zero temperature ($\beta \to \infty$) limit of \eqref{SrhoTr} gives
\begin{equation}\label{tildeTrlim}
\lim_{\beta \to \infty} \Tr^\oplus_{\beta,S}(A)
= \widetilde\Tr^\oplus(A) .
\end{equation}
\end{lem}

\proof We can identify $\cM^{(N)} = \cup_{U\in U(N)} U \cdot \Delta_{N-1}$, with the
simplex $\Delta_{N-1} =\{ P=(p_i)_{i=1}^N\,|\, p_i\geq 0, \, \sum_i p_i=1 \}$, where
the action of $U\in U(N)$ is by $P \mapsto U\cdot P:= U^* PU$, where $P$ is identified
with the diagonal density matrix with diagonal entries $p_i$. We then have
$$ \lim_{\beta \to \infty} \Tr^\oplus_{\beta,S}(A) = \min_{\rho \in \cM^{(N)}} \{ \Tr(\rho A) \} =
\min_{P=(p_i)\in \Delta_{N-1}}\{ \Tr(P \, UAU^*) \} $$ 
$$ = \min_{P=(p_i)\in \Delta_{N-1}} \{ \sum_i p_i (UAU^*)_{ii} \} = \min_i \{ (UAU^*)_{ii} \}. $$
\endproof

Recall that, for $\rho,\sigma\in \cM^{(N)}$, the quantum relative entropy is defined as
\begin{equation}\label{def:relS}
S(\rho || \sigma) = \Tr(\rho (\log\rho - \log\sigma)).
\end{equation}

\begin{lem}\label{relS}
When $A=A^*$ with $A\geq 0$, the expression $\Tr(\rho A) - \beta^{-1} \cN(\rho)$ can
identified with a relative entropy
\begin{equation}\label{eq:relS}
\Tr(\rho A) - \beta^{-1} \cN(\rho) = \frac{1}{\beta}  S(\rho || \sigma_{\beta,A}) - \frac{1}{\beta} \log Z_A(\beta),
\end{equation}
where
$$ \sigma_{\beta,A} = \frac{e^{-\beta A}}{Z_A(\beta)}, \ \  \ \text{ with } \ \ \  Z_A(\beta) =\Tr(e^{-\beta A}). $$
\end{lem}

\proof This follows by simply writing
$$ \Tr(\rho A )+ \beta^{-1} \Tr(\rho\log \rho) =\beta^{-1}  \Tr(\rho (\log\rho - \log e^{-\beta A})). $$
\endproof

\begin{prop}\label{ExpTr}
When $A=A^*$ with $A\geq 0$, the deformed trace \eqref{SrhoTr} with $S=\cN$ the von Neumann entropy 
is given by
\begin{equation}\label{ZdefTr}
\Tr^\oplus_{\beta,\cN}(A) = - \frac{\log Z_A(\beta)}{\beta},
\end{equation}
with $Z_A(\beta) =\Tr(e^{-\beta A})$.
\end{prop}

\proof By the previous lemma, we have 
$$ \Tr^\oplus_{\beta,\cN}(A) =\min_{\rho \in \cM^{(N)}} 
\frac{1}{\beta}  S(\rho || \sigma_{\beta,A}) - \frac{1}{\beta} \log Z_A(\beta) $$
The relative entropy have the property that $S(\rho ||\sigma)\geq 0$ with minimum
at $\rho =\sigma$ where $S(\rho || \rho)=0$. Thus, the minimum of the above expression
is $\beta^{-1} \log Z_A(\beta)$.
\endproof

\begin{rem} {\rm 
The expression $\beta^{-1} \log Z_A(\beta)$ can be interpreted as the Helmholtz
free energy in quantum statistical mechanics. }
\end{rem}

\begin{rem} {\rm
In the case where $A$ is the $2\times 2$ diagonal matrix with diagonal entries $(x,y)$,
the expression $\Tr^\oplus_{\beta,\cN}(A) =\beta^{-1} \log Z_A(\beta)$ recovers the
usual deformed addition $x\oplus_{\beta,S} y$ in the thermodynamic
semiring with $S$ the Shannon entropy, in the form \eqref{logexpbetaS},
$$ x\oplus_{\beta,S} y = \min_p \{ p x + (1-p)y - \beta^{-1} S(p) \} = - \beta^{-1} \log (e^{-\beta x} + e^{-\beta y}). $$ }
\end{rem}

\begin{cor}\label{Ulambda}
When $A=A^*$ with $A\geq 0$, the zero temperature limit is 
$$ \widetilde\Tr^\oplus(A)=\min\{ \lambda \in {\rm Spec}(A) \}. $$
\end{cor}

\proof We take the limit as $\beta \to \infty$ of $\Tr^\oplus_{\beta,\cN}(A) =\beta^{-1} \log Z_A(\beta)$.
The leading term is given by $\beta^{-1} \log e^{-\beta \lambda_{\rm min}}$, where 
$\lambda_{\rm min}=\min\{ \lambda \in {\rm Spec}(A) \}$, hence comparing with Lemma \ref{Trbeta0}, we get
$$ \min_{U\in U(N)} \min_i \{ (UAU^*)_{ii} \} = \lambda_{\rm min}. $$
\endproof

In the following, we will use the unconventional symbol $\boxplus$ for the direct sum of matrices,
to distinguish it from the symbol $\oplus$ that we have adopted for the addition operation
in min-plus semirings. The deformed trace has following behavior. 

\begin{prop}\label{tensTr}
For a matrix $A^*=A$, $A\geq 0$, that is a direct sum of two
matrices $A=A_1 \boxplus A_2$ with $A_i=A_i^*$ and $A_i\geq 0$, the deformed traces satisfy
\begin{equation}\label{odotprodTr}
\Tr^\oplus_{\beta,\cN}(A) = \Tr^\oplus_{\beta,\cN}(A_1) \odot \Tr^\oplus_{\beta,\cN}(A_2),
\end{equation}
where $\odot$ is the product in the tropical semiring $\bT$.
\end{prop}

\proof By Proposition \ref{ExpTr} we have $\Tr^\oplus_{\beta,\cN}(A) = - \beta^{-1} \log \Tr(e^{-\beta A})$.
For $A$ a direct sum of $A_1$ and $A_2$ we have $e^{-\beta A}=e^{-\beta A_1} \otimes e^{-\beta A_2}$
and $\Tr(e^{-\beta A_1} \otimes e^{-\beta A_2})=\Tr(e^{-\beta A_1}) \Tr(e^{-\beta A_2})$, hence we get
$$ \Tr^\oplus_{\beta,\cN}(A) = - \beta^{-1} \left(\log \Tr(e^{-\beta A_1}) + \log \Tr(e^{-\beta A_2}) \right). $$
\endproof

\medskip
\subsection{Generalizations of von Neumann entropy and relative entropy}

In addition to the von Neumann entropy, there are several other natural entropy functionals
in quantum information. Some of the main examples (see \cite{BenZyc}, \cite{Wilde}) are
\begin{itemize}
\item The quantum relative entropy: for $\rho,\sigma\in \cM^{(N)}$
$$ S(\rho || \sigma)= \Tr (\rho (\log\rho - \log\sigma)). $$
\item The quantum R\'enyi entropy: for $\rho \in \cM^{(N)}$
$$ S_q(\rho)=\frac{1}{1-q} \log \Tr(\rho^q). $$
\item The Belavkin--Staszewski relative entropy: for $\rho,\sigma\in \cM^{(N)}$
$$ S_{BS}(\rho || \sigma) = \Tr(\rho \log(\rho^{1/2} \sigma^{-1} \rho^{1/2})). $$
\item The quantum Tsallis entropy: for $\rho \in \cM^{(N)}$
$$ S_\alpha(\rho)= \frac{1}{1-\alpha}  \Tr(\rho(\rho^{\alpha-1}-I)). $$
\item The Umegaki deformed relative entropy: for $\rho,\sigma\in \cM^{(N)}$
$$ S_\alpha (\rho || \sigma)= \frac{4}{1-\alpha^2} \Tr((I-\sigma^{(\alpha+1)/2} \rho^{(\alpha-1)/2})\rho). $$
\end{itemize}
The fact that there is a large supply of entropies and relative entropies in the quantum case depends
on the fact that the expression $\rho \sigma^{-1}$, in going from classical to quantum, can be replaced 
by several different expressions, when $\rho$ and $\sigma$ do not commute. All of these entropy
functionals give rise to corresponding thermodynamic deformations  $\Tr^\oplus_{\beta, S}(A)$ 
of the tropical trace $\Tr^\oplus(A)$, defined as in  \eqref{SrhoTr}. In the case of relative entropies,
we assume given a fixed density matrix $\sigma$ and we set $S_\sigma(\rho)=S(\rho || \sigma)$, so that
$$ \Tr^\oplus_{\beta, S_\sigma}(A) = \min_{\rho \in \cM^{(N)}}  \{ \Tr(\rho A) -\beta^{-1} S(\rho || \sigma)\}. $$
This is the natural generalization of the case of thermodynamic semirings with $S$ 
the Kullback--Leibler relative entropy, discussed in \cite{MaThor}.

\medskip

\subsection{Thermodynamically deformed states on $C^*$-algebras}

The construction discussed above using entropy functionals on matrix algebras can
be extended to a more general setting of $C^*$-algebras. Let $\cA$ be a unital separable
$C^*$-algebra (noncommutative in general) and let $\cM$ be the convex set of states
on $\cA$, namely continuous linear functionals $\varphi: \cA \to \C$ that are normalized
by $\varphi(1)=1$ and satisfy the positivity condition $\varphi(a^* a) \geq 0$, for all $a\in \cA$.

\smallskip

There is a general notion of relative entropy $S(\varphi || \psi)$ of states $\varphi,\psi\in \cM$ 
on a $C^*$-algebra $\cA$, see \S 2.3 of \cite{NeStor}, with the same bi-convexity property of the
usual relative entropy, and with $S(\varphi || \psi)\geq 0$ for all states $\varphi, \psi$, with 
equality attained only when $\varphi=\psi$.
In the case where there is a trace $\tau: \cA \to \C$ on the $C^*$-algebra, consider
states of the form $\varphi(a)=\tau(a \xi)$, $\psi(a)=\tau(a \eta)$, where $\xi,\eta$ are positive
elements in the algebra with $\tau(\xi)=\tau(\eta)=1$. Then the relative entropy reduces to 
\begin{equation}\label{relStau}
S(\varphi || \psi) = \tau (\xi (\log \xi - \log \eta)).
\end{equation}

Given a state $\psi \in \cM$, we define its thermodynamical deformation $\psi_{\beta,S}$ as
\begin{equation}\label{psibetaS}
\psi_{\beta,S}(a) = \min_{\varphi \in \cM} \{ \varphi (a) + \beta^{-1} S(\varphi || \psi) \}. 
\end{equation}

Notice that, in the finite dimensional case of a matrix algebra, this agrees with our
previous definition of the deformation of the trace, since states are of the form $\varphi(a)=\Tr(a \rho)$
for some density matrix $\rho$ and the von Neumann entropy can be seen as
$\cN(\rho)=-S(\varphi || \psi)$, for $\varphi(a)=\Tr(a \rho)$ and $\psi(a)=\Tr(a)$.

\medskip

While this more general setting is not the main focus of the present paper, we illustrate the 
construction in one significant example. Let $\cA_\theta$ be the irrational rotation algebra (noncommutative
torus) with unitary generators $U,V$ satisfying $VU=e^{2\pi i \theta} UV$. Let 
$\tau$ be the canonical trace, $\tau(U^n V^m)=0$ for $(n,m)\neq (0,0)$ and $\tau(1)=1$.
We consider only states of the form $\varphi(a)=\tau(a \xi)$ for some positive element
$\xi \in \cA_\theta$. Let $\cM_\tau$ be the set of such states. We then consider the thermodynamic
deformation of the canonical trace given by
\begin{equation}\label{taubetaS}
\tau_{\beta,S}(a) = \min_{\varphi \in \cM_\tau} \{ \varphi(a) + \beta^{-1} S(\varphi || \tau) \}.
\end{equation}

We then obtain the following result, whose proof is completely analogous to Proposition \ref{ExpTr} above.

\begin{prop}\label{NCtorus}
For $a=h^*h \geq 0$ in $\cA_\theta$, the deformed trace $\tau_{\beta,S}(a)$ is given by
$$ \tau_{\beta,S}(a) = \min_{\varphi \in \cM_\tau} \{ \beta^{-1} S(\varphi || \varphi_{\beta,a}) -\beta^{-1} \log \tau (e^{-\beta a}) \}=
-\beta^{-1} \log \tau (e^{-\beta a}) $$
for the  KMS$_\beta$ state
$$ \varphi_{\beta,a}(b) =\frac{\tau(b e^{-\beta a})}{\tau(e^{\beta a})} $$
of the time evolution $\sigma_t(b) =e^{it a} b e^{-ita}$ on $\cA_\theta$,
with $-\beta^{-1} \log \tau (e^{-\beta a})$ the associated 
the Helmholtz free energy.
\end{prop}

The limit $\lim_{\beta\to \infty} \tau_{\beta,S}(a)$ should then be regarded as a notion of 
``tropicalization" of the von Neumann trace $\tau$ of the noncommutative torus.

\medskip
\subsection{Thermodynamic deformations and entropy of dynamical systems}

We consider a locally compact Hausdorff space $X$, with a dynamical
system $\sigma: X \to X$. We focus in particular on the case
where $X$ is a Cantor set, identified with the set of infinite words
$w =w_0 w_1 \ldots w_i w_{i+1} \ldots$ in a finite alphabet $w_i \in \fA$,
with $\# \fA=n$, with the topology generated by cylinder sets 
$\cC(a_0,\ldots, a_N)=\{ w \in X\,|\, w_i=a_i, \, 0\leq i \leq N \}$. Let $d(x,y)$
be a compatible metric. As dynamical system, we consider in particular
the case of the one-sided shift $\sigma: X \to X$, defined by $\sigma(w)_i = w_{i+1}$.

\smallskip

A Bernoulli measure $\mu_P$ on $X$ is a shift-invariant measure defined by a probability
$P=(p_1,\ldots, p_n)$, with $p_i\geq 0$ and $\sum_{i=1}^n p_i=1$, on the alphabet $\fA$.
It assigns measure $\mu_P(\cC(a_0,\ldots, a_N))=p_{a_0}\cdots p_{a_N}$ to the cylinder sets.

\smallskip

A Markov measure $\mu_{P,\rho}$ on $X$ is a shift-invariant measure defined by a pair $(P,\rho)$
of a probability $P=(p_1,\ldots, p_n)$ on $\fA$ and a stochastic matrix $\rho$ satisfying $P \rho=P$.
It assigns measure $\mu_{P,\rho}(\cC(a_0,\ldots, a_N))= p_{a_0} \rho_{a_0a_1}\cdots \rho_{a_{N-1} a_N}$.
A Markov measure $\mu_{P,\rho}$ is supported on a subshift of finite type $X_A \subset X$,
given by $X_A=\{ w \in X \,|\, A_{w_i w_{i+1}}=1, \, \forall i\geq 0 \}$, where the matrix $A_{ij}$
has entries $0$ or $1$, according to whether the corresponding entry $\rho_{ij}$ of the
stochastic matrix $\rho$ is $\rho_{ij}=0$ or $\rho_{ij}\neq 0$. The subspace $X_A$ is shift-invariant.

\smallskip

Recall that, for $\mu$ a $\sigma$-invariant probability measure on $X$, one defines the entropy 
$S(\mu,\sigma)$ as the $\mu$-almost everywhere value of the local entropy
$$ h_{\mu,\sigma}(x) = \lim_{\delta \to 0} \lim_{n \to \infty} - \frac{1}{n} \log \mu(B_\sigma (x,n,\delta)), $$
where $B_\sigma(x,n,\delta)=\{ y \in X\,|\, d(\sigma^j(x), \sigma^j(y))< \delta, \, \forall 0\leq j \leq n \}$
are the Bowen balls, see \cite{PeCli}. In the case of a Bernoulli measure $\mu=\mu_P$, the dynamical
entropy agrees with the Shannon entropy of $P$,
$$ S(\mu_P,\sigma)=-\sum_{i=1}^N p_i \log p_i , $$ 
while for a Markov measure, the dynamical entropy is
$$ S(\mu_{P,\rho},\sigma)=-\sum_{i=1}^N p_i \sum_{j=1}^N \rho_{ij} \log \rho_{ij} . $$

\smallskip

In the same spirit as the thermodynamic deformations of the trace discussed 
previously in this section, we can introduce thermodynamic deformations of
the integral of functions $f\in C(X,\R)$ by setting
\begin{equation}\label{intbetaS}
\int_X^{(\beta,S)} f(x) dx := \inf_\mu \{ \int_X f(x) d\mu(x) - \beta^{-1} S(\mu,\sigma) \},
\end{equation}
where the infimum is taken over a specific class of $\sigma$-invariant measures, for
example over all Bernoulli measures, or over all Markov measures, or more generally over
all $\sigma$-invariant ergodic measures. In the latter case, recall that the topological
entropy of the shift $\sigma$ is
$$ h(X,\sigma) = \sup_\mu \{ S(\mu,\sigma) \}, $$
with the supremum taken over all $\sigma$-invariant ergodic measures. 

\smallskip

We will not discuss further the properties of the functionals \eqref{intbetaS}, as that
would lead us outside the main scope of the present paper.

\bigskip
\section{Rota--Baxter structures and Birkhoff factorization in min-plus semirings}\label{RBBirkSec}

A mathematical model of renormalization for perturbative quantum field
theories, based on a commutative Hopf algebra $\cH$, a Rota--Baxter algebra $\cR$ 
and the Birkhoff factorization of morphisms of commutative algebras 
$\phi: \cH \to \cR$, was developed in \cite{CoKr}, \cite{EFGuo}, \cite{EFGK}.
More recently, an approach to the theory of computation and the halting
problem modelled on quantum field theory and renormalization was
developed in \cite{Man1}, \cite{Man2}, \cite{Man3}, and further investigated in
\cite{DelMar}. In view of applications to the theory of computation, it was observed
in \S 4.6 of \cite{Man1} that it would be useful to replace characters given by
commutative algebra homomorphisms $\phi: \cH \to \cR$ from the Hopf algebra
to a Rota--Baxter algebra, with characters $\psi: \cH \to \bS$ with values in a
min-plus semiring, satisfying $\psi(xy)=\psi(x)+\psi(y)=\psi(x)\odot \psi(y)$. With
this motivation in mind, we develop here a setting for Rota--Baxter structures 
and Birkhoff factorization taking place in min-plus semirings and in their
thermodynamic deformations.

\medskip
\subsection{Rota--Baxter algebras and renormalization}

We refer the reader to \cite{Guo} for a general introduction to the
subject of Rota--Baxter algebras. For their use in 
renormalization of perturbative quantum field theories, we refer
the reader to \cite{CoMa}, \cite{EFGK}, \cite{Mar}, for more details.

\smallskip

A Rota--Baxter algebra (ring) of weight $\lambda$ is a unital commutative algebra
(ring) $\cR$ endowed with a linear operator $\cT:\cR \to \cR$ which satisfies the 
$\lambda$-Rota--Baxter identity
\begin{equation}\label{RBlambda}
 \cT(a) \cT(b) = \cT(a \cT(b))+ \cT(\cT(a)b) +\lambda \cT(ab) .
\end{equation} 

\smallskip

We will be especially interested in two cases, namely $\lambda=\pm 1$, which
correspond, respectively, to the identities
\begin{equation}\label{RBplus1}
\cT(a) \cT(b) = \cT(a \cT(b)) + \cT(\cT(a) b) + \cT(ab) ,
\end{equation}
\begin{equation}\label{RBminus1}
\cT(a) \cT(b)  + \cT(ab) = \cT(a \cT(b)) + \cT(\cT(a) b).
\end{equation}
The latter case, with weight $\lambda=-1$, is the one used in renormalization
in quantum field theory, while we will see that the case $\lambda=+1$ is more
natural to adapt to the setting of min-plus semirings.

\smallskip

Laurent polynomials $\cR=\C[t,t^{-1}]$ with the projection $\cT$
onto the polar part are the prototype example of a Rota--Baxter 
algebra of weight $-1$. Recall also that, if $\cT$ is a Rota--Baxter
operator of weight $\lambda\neq 0$, then $\lambda^{-1} \cT$ is a 
Rota--Baxter operator of weight $1$.

\smallskip

When $\lambda=-1$ the Rota--Baxter operator $\cT$ determines
a decomposition of $\cR$ into two commutative algebras (rings),
$\cR_+=(1-\cT)\cR$ and $\cR_-$ given by the unitization of $\cT\cR$.

\smallskip

Algebraic renormalization is a factorization procedure for
Hopf algebra characters. More precisely, one considers over
a field or ring $k$ a graded connected commutative Hopf algebra 
$\cH =\oplus_{n\geq 0} \cH_n$ with $\cH_0=k$ and the
set of homomorphisms of commutative rings $\Hom(\cH,\cR)$,
where the target $\cR$ is a Rota--Baxter ring of weight $\lambda=-1$.

\smallskip

The convolution product $\star$ of morphisms $\phi_1, \phi_2 \in \Hom(\cH,\cR)$
is dual to the coproduct in $\cH$, that is, 
\begin{equation}\label{phiconvol}
\phi_1\star \phi_2(x) = \langle \phi_1\otimes \phi_2, \Delta(x)\rangle
=\sum \phi_1(x^{(1)}) \phi_2(x^{(2)}),
\end{equation}
where 
$$ \Delta(x)=\sum x^{(1)}\otimes x^{(2)}=x \otimes 1 + 1 \otimes x + \sum x' \otimes x''. $$

\smallskip

The Birkhoff factorization of a morphism $\phi\in \Hom(\cH,\cR)$ is
a multiplicative decomposition
\begin{equation}\label{Birkhoff}
\phi = (\phi_-\circ S)\star \phi_+ ,
\end{equation}
where  $S$ is the antipode, defined inductively by
$$ S(x)=-x -\sum S(x') x'', $$
where $\Delta(x)=x\otimes 1 + 1 \otimes x + \sum x' \otimes x''$,
with the $x'$, $x''$ of lower degrees. 

\smallskip

The two parts $\phi_\pm$ of the
Birkhoff factorization are morphisms of commutative algebras (rings)
$\phi_\pm: \cH \to \cR_\pm$. The decomposition is obtained inductively
through the explicit formula
\begin{equation}\label{recBirkhoff}
 \phi_-(x)=-\cT (\phi(x) +\sum \phi_-(x') \phi(x''))  \ \ \text{ and } \ \ 
  \phi_+(x)=(1-\cT)(\phi(x) +\sum \phi_-(x') \phi(x'')) .
\end{equation}  
We denote by $\tilde\phi(X)$ the Bogolyubov-Parashchuk ``preparation" of $\phi(X)$
\begin{equation}\label{BPprep}
 \tilde\phi(x) :=\phi(x) +\sum \phi_-(x') \phi(x'') . 
\end{equation} 

\smallskip

The fact that $\phi_-$, constructed inductively as above, 
is still a homomorphism of commutative rings  $\phi_-: \cH \to \cR_-$ is obtained
by comparing
\begin{equation}\label{phiminXY}
 \phi_-(xy)=- \cT(\tilde\phi(x) \tilde\phi(y)) + \cT( \cT(\tilde\phi(x)) \tilde\phi(y) +
\tilde\phi(x) \cT(\tilde\phi(y)) 
\end{equation}
and
\begin{equation}\label{phiminXY2}
 \phi_-(x)\phi_-(y) = \cT(\tilde\phi(x)) \cT(\tilde\phi(y)) 
\end{equation}
using the Rota--Baxter identity for $\cT$. It then follows that $\phi_+$ is
also a ring homomorphism. The expression for $\phi_-(xy)$ above
is easily obtained by decomposing the terms $(xy)'$ and $(xy)''$ in the
non-primitive part of the coproduct $\Delta(xy)$ in terms of $x$, $y$,
$x'$ and $x''$, $y'$ and $y''$. We will return to this argument below.

\medskip
\subsection{Rings and semirings}

The usual setting of Rota--Baxter algebras recalled above is based on
commutative rings $\cR$ with a linear operator $\cT$ satisfying the
identity \eqref{RBlambda}. Our purpose in this section is to extend
this notion to min-plus semirings and their thermodynamic deformations
and relate the Rota--Baxter structures and Birkhoff factorization on
semirings to the ordinary ones on rings.

\smallskip

To this purpose, we will consider min-plus semirings, and thermodynamic
deformations that are related to commutative ring via a ``logarithm" map.

\smallskip

The kind of semirings we consider are semirings $\bS$ 
with min-plus operations $\oplus$, $\otimes$, for which
thermodynamic deformations $\bS_{\beta,S}$ are defined,
with $S$ the Shannon entropy.

\smallskip

The Gelfand correspondence between compact Hausdorff spaces $X$ and
commutative unital $C^*$-algebras $C(X)$ admits a generalizations for
semirings of continuous functions $C(X,\bT)$ with values in the tropical semiring
$\bT=(\bR\cup \{ \infty \}, \oplus, \odot)$, with the pointwise operations, see \cite{KerSchn}. 

\smallskip

Thus, a large class of examples of semirings $\bS$ of the type described
above is given by $\bS=C(X,\bT)$ with the pointwise $\oplus=\min$ and
$\odot=+$ operations, and their thermodynamic deformations $\bS_{\beta,S} =C(X, \bT_{\beta,S})$
with the pointwise deformed addition $\oplus_{\beta,S}$.

\smallskip

\begin{defn}\label{RandS}
Let $\cR$ be a commutative ring (or algebra) and let $\bS$ be a semiring
with min-plus operations $\oplus$, $\otimes$. The pair $(\cR, \bS)$ is a {\em logarithmically related}
pair, if there is a bijective map $\cL: {\rm Dom}(\cL)\subset \cR \to \bS$ satisfying 
$\cL(ab)=\cL(a)+\cL(b)=\cL(a) \odot \cL(b)$, for all $a,b \in {\rm Dom}(\cL)$.
\end{defn}

\smallskip

Let $\bS_{\beta,S}$ be a thermodynamic deformation of $\bS$, 
for which we can write the deformed addition as
$$ f_1 \oplus_{\beta,S} f_2 = -\beta^{-1} \log ( E(-\beta f_1) + E(-\beta f_2) ), $$
where $E: \bS \to {\rm Dom}(\cL)\subset \cR$ denotes the inverse of $\cL$ and
$+$ is addition in the ring $\cR$. The undeformed $\odot$ operation
in $\bS$ is related to the product in $\cR$ by $f_1 + f_1 = \cL(E(f_1))+\cL(E(f_2))=\cL(E(f_1) E(f_2))$.

\smallskip

\begin{ex}{\rm 
The above applies to the case for all the semirings $\bS_{\beta,S}=C(X,\bT_{\beta,S})$, with the usual
deformed addition, with $S$ the Shannon entropy, given by 
$$ f_1 \oplus_{\beta,S} f_2 = -\beta^{-1} \log ( e^{-\beta f_1} + e^{-\beta f_2} ). $$
In this case $\cR= C(X,\R)$ and the subset ${\rm Dom}(\cL)\subset \cR$ is given by functions
$a \in C(X,\R^*_+)$, with $a=e^{-\beta f}$. In other words $\cL(a) = -\beta^{-1} \log(a)$.}
\end{ex}

\smallskip

\begin{ex}\label{formlog}{\rm
Let $\cR$ be the ring of formal power series $\Q[[t]]$ in one variable, with rational coefficients.
Let ${\rm Dom}(\cL)\subset \cR$ be the subset of power series $\alpha(t)=\sum_{k\geq 0} a_k t^k$ 
with $a_0=1$. Then $\cL$ is the formal logarithm
$\cL(1+\alpha) =\alpha - \frac{1}{2} \alpha^2 + \frac{1}{3} \alpha^3 + \cdots =  
\sum_{k=1}^\infty \frac{(-1)^{k+1}}{k} \alpha^k$, mapping $\cL: {\rm Dom}(\cL) \to \Q[[t]]$.
It satisfies $\cL(\alpha \gamma)=\cL(\alpha) + \cL(\gamma)$. The inverse of the formal logarithm $\cL$
is given by the formal exponential $E(\gamma)=\sum_{k\geq 0} \gamma^k/k!$. We can view $\Q[[t]]$
as a thermodynamic semiring with deformed addition 
$$ \alpha_1 \oplus_{\beta,S} \alpha_2 = \beta^{-1} \cL ( E(-\beta \alpha_1) + E (-\beta \alpha_2)), $$
for $\beta\in \Q$, and undeformed multiplication $\alpha_1 \odot \alpha_2 = \alpha_1 + \alpha_2$.}
\end{ex}

\smallskip

For simplicity, in the following we will always write simply $\log$ and $\exp$ for the maps
relating a ring $\cR$ and a semiring $\bS$, as in definition \ref{RandS}.

\medskip
\subsection{Birkhoff factorization in min-plus semirings}

\smallskip

\begin{defn}\label{charminplus}
Let $\bS$ be a min-plus semiring as above. Let $\cH$ be a graded, connected, commutative
Hopf algebra. A {\em min-plus character} (or {\em $\bS$-character}) of the Hopf algebra is a map 
$\psi: \cH \to \bS$ that satisfies the conditions  $\psi(1)=0$ and
\begin{equation}\label{psiadd}
\psi(xy)= \psi(x) + \psi(y), \ \ \  \forall x,y \in \cH .
\end{equation}
\end{defn}

We also define a convolution product of min-plus characters. The intuition
behind the definition comes from a standard heuristic 
reasoning, which regards the min-plus algebra as the ``arithmetic of orders
of magnitude". Namely, when $\epsilon\to 0$, the leading term in 
$\epsilon^\alpha + \epsilon^\beta$ is $\epsilon^{\min\{\alpha,\beta\}}$,
while the leading term of $\epsilon^\alpha \epsilon^\beta$ is
$\epsilon^{\alpha+\beta}$. Thus, the notion of convolution product for
min-plus characters should reflect the behavior of the leading order in
the usual notion of convolution product of (commutative algebra valued)
characters. 

\begin{defn}\label{convolminplus}
For $\psi_1, \psi_2$ as above, the convolution product $\psi_1\star \psi_2$
is given by
\begin{equation}\label{convolpsi}
(\psi_1 \star \psi_2)(x)=\min\{ \psi_1(x^{(1)}) + \psi_2(x^{(2)}) \} = \bigoplus ( \psi_1(x^{(1)}) \odot \psi_2(x^{(2)}) ),
\end{equation}
where the minimum is taken over all the pairs $(x^{(1)}, x^{(2)})$ that appear in
the coproduct $\Delta(x)=\sum x^{(1)}\otimes x^{(2)}$ in the Hopf algebra $\cH$,
and $\oplus=\min$ and $\odot=+$ are the (pointwise) operations of the semiring $\bS$.
\end{defn}

Similarly, we reformulate the notion of Birkhoff factorization in the following way.

\begin{defn}\label{semiBirk}
Let $\psi$ be a min-plus character of the Hopf algebra $\cH$. A Birkhoff factorization
of $\psi$ is a decomposition $\psi_+=\psi_- \star \psi$, with $\star$ the convolution product \eqref{convolpsi},
where $\psi_\pm$ satisfy \eqref{psiadd}.
\end{defn}

Notice that, unlike the usual way of writing Birkhoff factorizations in the form \eqref{Birkhoff},
the formulation above as $\psi_+=\psi_- \star \psi$ does not require the use of the antipode of
the Hopf algebra, hence it extends to the case where $\cH$ is a bialgebra. Since in our main
applications $\cH$ will be a Hopf algebra, we maintain this assumption in the following.

\medskip
\subsection{Rota--Baxter operators on min-plus semirings}

Let $\bS$ be a min-plus semiring, with (pointwise) operations $\oplus$ and $\odot$.
A map $T: \bS \to \bS$ is $\oplus$-additive if it is
monotone, namely $T(a)\leq T(b)$  if $a\leq b$, for all $a,b \in \bT$.
For a semiring of the form $\bS=C(X,\bT)$ the condition is 
pointwise in $t\in X$.

\smallskip

We define Rota--Baxter structures with weight $\lambda>0$ as follows.

\smallskip

\begin{defn}\label{RBsemiplus}
A min-plus semiring $(\bS, \oplus, \odot)$ is a Rota--Baxter semiring of
weight $\lambda>0$ if there is a $\oplus$-additive map $T: \bS \to \bS$, which 
for all $f_1,f_2 \in \bS$ satisfies the identity
\begin{equation}\label{RBsemieqplus}
  T(f_1) \odot T(f_2)  = T(T(f_1) \odot f_2) \oplus T(f_1 \odot T(f_2)) \oplus T(f_1 \odot f_2) \odot \log \lambda .
\end{equation} 
\end{defn}

Similarly, we can define Rota--Baxter structures of weight $\lambda <0$ in the following way.

\begin{defn}\label{RBsemimin}
A min-plus semiring $(\bS, \oplus, \odot)$ is a Rota--Baxter semiring of
weight $\lambda <0$ if there is a $\oplus$-additive map $T: \bS \to \bS$, which 
for all $f_1,f_2 \in \bS$ satisfies the identity
\begin{equation}\label{RBsemieqmin}
  T(f_1) \odot T(f_2) \oplus T(f_1 \odot f_2) \odot \log (-\lambda) = T(T(f_1) \odot f_2) \oplus T(f_1 \odot T(f_2)).
\end{equation} 
\end{defn}

\smallskip

We have the following result on the existence of Birkhoff factorizations.
As in the usual case, the proof is constructive, as it inductively defines
the two parts of the factorization.

\begin{thm}\label{algrensemi}
Let $\psi: \cH \to \bS$ be a min-plus character of a 
graded, connected, commutative Hopf algebra $\cH$.
Assume that the target semiring $\bS$ has a Rota--Baxter structure of weight $+1$,
as in Definition \ref{RBsemiplus}. Then there is a Birkhoff factorization $\psi_+ = \psi_- \star \psi$,
where $\psi_-$ and $\psi_+$ are also min-plus characters.
\end{thm}

\proof As in the usual Rota--Baxter algebra case, we construct the factors
$\psi_\pm$ inductively.  We define the Bogolyubov-Parashchuk preparation
of $\psi$ as
\begin{equation}\label{BPpsi}
\tilde\psi(x)=\min \{ \psi(x), \psi_-(x') + \psi(x'') \} = \psi(x) \oplus \bigoplus \psi_-(x') \odot \psi(x''),
\end{equation}
where $(x',x'')$ ranges over all pairs in the non-primitive
part of the coproduct $\Delta(x) =x\otimes 1 + 1 \otimes x + \sum x' \otimes x''$,
and $\psi_-$ is assumed defined by induction on all the lower degree terms 
$x'$ in the Hopf algebra as
\begin{equation}\label{psiminus}
\psi_-(x) := T(\tilde\psi(x)) = T( \min\{ \psi(x), \psi_-(x') + \psi(x'') \} ) = T\left( \psi(x) \oplus \bigoplus \psi_-(x') \odot \psi(x'') \right). 
\end{equation}
By the $\oplus$-linearity of $T$, this is the same as
$$ \psi_-(x) =  \min \{ T(\psi(x)), 
T( \psi_-(x') + \psi(x'') ) \} = T(\psi(x)) \oplus \bigoplus T( \psi_-(x') \odot \psi(x'') ) . $$
The positive part of the factorization is then obtained as the convolution product
\begin{equation}\label{psiplus}
\psi_+(x) :=(\psi_- \star \psi)(x)=\min\{ \psi_-(x), \psi(x), \psi_-(x')+\psi(x'') \} =
\min\{ \psi_-(x), \tilde\psi(x) \} = \psi_-(x) \oplus \tilde\psi(x).
\end{equation}
We need to check that $\psi_\pm$ satisfy \eqref{psiadd}.
We have $\psi_-(xy)= T \min\{ \psi(x)+\psi(y), \psi_-((xy)') + \psi((xy)'') \}$, where
we can decompose the terms $(xy)'$ and $(xy)''$ in terms of
$x$, $y$, $x'$ and $x''$, $y'$ and $y''$. This gives 
\begin{equation}\label{psiminxylist}
\psi_-(xy)= T \min\left\{ \begin{array}{l} \psi(x)+\psi(y), \\ 
\psi_-(x)+\psi(y), \\ \psi_-(y)+\psi(x), \\ \psi_-(y')+\psi(x y'') , \\
\psi_-(x')+\psi(x'' y) , \\ \psi_-(x y')+ \psi(y''), \\ \psi_-(x' y)+\psi(x''), \\
\psi_-(x'y') + \psi(x''y'')\end{array} \right\} . 
\end{equation}
Using associativity and commutativity of $\oplus$ and $\oplus$-additivity of $T$,
we can group these terms together into
$$ \psi_-(xy)= \min \{ \alpha(x,y,x',y'), \beta(x,y,x',y') \}, $$
where we have 
\begin{equation}\label{min1part}
\alpha(x,y,x',y') =T \min \{ \psi_-(x)+\psi(y), \psi(x) +\psi_-(y),
\psi_-(xy')+\psi(y''), \psi_-(x'y)+\psi(x'') \}
\end{equation}
\begin{equation}\label{min2part}
\beta(x,y,x',y') = T \min\{ \psi(x)+\psi(y), \psi_-(y')+\psi(x y''), \psi_-(x')+\psi(x''y), \psi_-(x'y')+\psi(x''y'') \}
\end{equation}
Assuming inductively that 
$$ \psi_-(uv) = \psi_-(u)+\psi_-(v), $$
for all terms $u$ and $v$ in $\cH$ of degrees 
$\deg(u)+\deg(v)< \deg(xy)$, 
and using the fact that $T$ is $\oplus$-additive,
we can rewrite the term $\alpha(x,y,x',y')$ of \eqref{min1part} as
\begin{eqnarray}\label{min1part1}
\alpha(x,y,x',y') & = & T \min\{ \psi_-(x)+\tilde\psi(y), \tilde\psi(x)+\psi_-(y) \} \\[2mm] \nonumber
& = & \min \{ T(T(\tilde\psi(x))+
\tilde\psi(y)), T(\tilde\psi(x)+T(\tilde\psi(y))) \}
\end{eqnarray}
and we can write the term $\beta(x,y,x',y')$ of \eqref{min2part} as
\begin{equation}\label{min2part2}
\beta(x,y,x',y') = T \min\{ \tilde\psi(x) + \tilde\psi(y) \} = \min \{ T(\tilde\psi(x) + \tilde\psi(y)) \}.
\end{equation}
Thus, we have
\begin{eqnarray}\label{psiminxyRB}
\psi_-(xy) &= & \min \{ T(\tilde\psi(x) + \tilde\psi(y)), T(T(\tilde\psi(x))+
\tilde\psi(y)), T(\tilde\psi(x)+T(\tilde\psi(y))) \} \\[2mm] \nonumber
& = & T(\tilde\psi(x) \odot \tilde\psi(y)) \oplus
T(T(\tilde\psi(x))\odot \tilde\psi(y)) \oplus T(\tilde\psi(x)\odot T(\tilde\psi(y))).
\end{eqnarray}
Since the operator $T$ satisfies the Rota--Baxter identity \eqref{RBsemieqplus}
with $\lambda=1$, we can rewrite the above as
$$ \psi_-(xy) = T(\tilde\psi(x)) \odot T(\tilde\psi(y)) = T(\tilde\psi(x)) + T(\tilde\psi(y)) = \psi_-(x) + \psi_-(y). $$
The fact that $\psi_+(xy)=\psi_+(x) + \psi_+(y)$ then follows from $\psi_+ = \psi_- \star \psi$.
\endproof

\bigskip

\section{Thermodynamic Rota--Baxter structures and Birkhoff factorizations}\label{thermoRBsec}

In Theorem \ref{algrensemi} we have used the associativity and commutativity 
properties of the tropical addition $\oplus$, in reordering the terms in $\psi_-(xy)$ to
prove it satisfies $\psi_-(xy)=\psi_-(x) + \psi_-(y)$.  Thus, in extending the result
to thermodynamic semirings, we will focus on the case of 
thermodynamic deformations $\oplus_{\beta,S}$, where 
$S$ is the Shannon entropy, since in this case both  
associativity and commutativity continue to hold for the deformed addition 
$\oplus_{\beta,S}$. 

\smallskip

\begin{defn}\label{defbetalin}
Let  $\bS_{\beta,S}$ be thermodynamic deformations  of a
semiring $\bS$, with operations
$\oplus_{\beta,S}$ and $\odot$, and with $S$ the Shannon entropy.
An operator $T: \bS_{\beta,S} \to \bS_{\beta,S}$ is $\oplus_{\beta,S}$-linear if, for all $f_1,f_2 \in \bS_{\beta,S}$
and all $\alpha,\gamma \in \bT$,
\begin{equation}\label{oplusbetaSlinear}
T( \alpha \odot f_1 \oplus_{\beta,S} \gamma \odot f_2)= \alpha \odot T(f_1)  \oplus_{\beta,S} \gamma \odot T(f_2).
\end{equation}
\end{defn}

\medskip
\subsection{Classical and thermodynamic Rota--Baxter operators}

As in the case of min-plus semirings $\bS=C(X,\bT)$ with pointwise $\oplus$ and $\odot$
operations, we can similarly define Rota--Baxter structures on their thermodynamic
deformations $\bS_{\beta,S}$. In the case of weight $\lambda >0$ we have the following.

\begin{defn}\label{RBplusbetaS}
A thermodynamic semiring $\bS_{\beta,S}$ is a Rota--Baxter semiring of
weight $\lambda>0$ if there is a $\oplus_{\beta,S}$-additive map $T: \bS_{\beta,S} \to \bS_{\beta,S}$, which 
for all $f_1,f_2 \in \bS_{\beta,S}$ satisfies the identity
\begin{equation}\label{RBplusbeta}
  T(f_1) \odot T(f_2)  = T(T(f_1) \odot f_2) \oplus_{\beta,S} T(f_1 \odot T(f_2)) \oplus_{\beta,S} T(f_1 \odot f_2) \odot 
  \log \lambda .
\end{equation} 
\end{defn}

The case with $\lambda <0$ is analogous: we have the following.

\begin{defn}\label{RBminbetaS}
A thermodynamic semiring $\bS_{\beta,S}$ is a Rota--Baxter semiring of
weight $\lambda <0$ if there is a $\oplus_{\beta,S}$-additive map $T: \bS_{\beta,S} \to \bS_{\beta,S}$, which 
for all $f_1,f_2 \in \bS_{\beta,S}$ satisfies the identity
\begin{equation}\label{RBsemieqmin}
  T(f_1) \odot T(f_2) \oplus_{\beta,S} T(f_1 \odot f_2) \odot  \log (-\lambda) = T(T(f_1) \odot f_2) \oplus_{\beta,S} 
  T(f_1 \odot T(f_2)).
\end{equation} 
\end{defn}

\smallskip

\begin{thm}\label{RBexp}
Let $\cR$ be a commutative ring and $\bS$ a min-plus semiring, logaritmically related as in Definition \ref{RandS}. 
Given $T : \bS\to \bS$, define a new map $\cT: \cR \to \cR$ by setting
$$\cT(e^{-\beta f}):= e^{-\beta T(f)},$$
for $a=e^{-\beta f}$ in ${\rm Dom}(\log)\subset \cR$. 
Then $T$ satisfies the Rota--Baxter identity \eqref{RBplusbeta} or \eqref{RBsemieqmin} of weight $\lambda$
if and only if $\cT$ satisfies the ordinary Rota--Baxter identity 
$$\cT(e^{-\beta f_1}) \cT(e^{-\beta f_2}) =\cT(\cT(e^{-\beta f_1}) e^{-\beta f_2}) + 
\cT( e^{-\beta f_1} \cT(e^{-\beta f_2}))+ \lambda_\beta\, \cT(e^{-\beta f_1}e^{-\beta f_2}).$$
of weight $\lambda_\beta=\lambda^{-\beta}$, for $\lambda>0$, or weight $\lambda_\beta=- |\lambda|^{-\beta}$ for
$\lambda<0$.
\end{thm}

\proof In the case $\lambda>0$, we write the left-hand-side of the Rota--Baxter identity of weight $\lambda_\beta$ for $\cT$ as
$$ \cT(e^{-\beta f_1}) \cT(e^{-\beta f_2}) = e^{-\beta (T(f_1)+T(f_2))} $$
while the right-hand-side gives
$$ \cT(\cT(e^{-\beta f_1)} e^{-\beta f_2}) + \cT (e^{-\beta f_1}) e^{-\beta f_2}) + \lambda_\beta \cT(e^{-\beta(f_1+f_1)}) $$
$$ =\cT(e^{-\beta (T(f_1)+f_2)} )+\cT(e^{-\beta (f_1+T(f_2))}) 
+ \lambda_\beta e^{-\beta T(f_1+f_2)} $$
$$ = e^{-\beta T(T(f_1)+f_2)}+ e^{-\beta T(f_1+T(f_2))} + 
e^{-\beta (T(f_1+f_2) -\beta^{-1} \log \lambda_\beta) } . $$
This gives the identity
$$ T(f_1)+T(f_2) = -\beta^{-1} \log( e^{-\beta T(T(f_1)+f_2)}+ e^{-\beta T(f_1+T(f_2))} + 
e^{-\beta (T(f_1+f_2) -\beta^{-1} \log \lambda_\beta) }). $$
For $\lambda_\beta=\lambda^{-\beta}$, this is equivalently written as
$$ T(f_1) \odot T(f_2)  = T(T(f_1) \odot f_2) \oplus_{\beta,S} T(f_1 \odot T(f_2)) \oplus_{\beta,S} T(f_1 \odot f_2) \odot 
\log \lambda. $$
In the case with $\lambda<0$, we write the left-hand-side of the Rota--Baxter identity for $\cT$ as
$$ \cT(e^{-\beta f_1}) \cT(e^{-\beta f_2}) - \lambda_\beta \cT(e^{-\beta f_1}e^{-\beta f_2}) =
e^{-\beta (T(f_1)+T(f_2))} + e^{-\beta (T(f_1+f_2) -\beta^{-1} \log (-\lambda_\beta)) }  $$
and the right-hand-side
$$ \cT(e^{-\beta (T(f_1)+f_2)} )+\cT(e^{-\beta (f_1+T(f_2))}) = e^{-\beta T(T(f_1)+f_2)}+ e^{-\beta T(f_1+T(f_2))} . $$
This gives the identity
$$ -\beta^{-1} \log(e^{-\beta (T(f_1)+T(f_2))} + e^{-\beta (T(f_1+f_2) -\beta^{-1} \log (-\lambda_\beta)) } ) =
-\beta^{-1} \log( e^{-\beta T(T(f_1)+f_2)}+ e^{-\beta T(f_1+T(f_2))} ). $$
For $\lambda_\beta=- |\lambda|^{-\beta}$, this is equivalently written as
$$ T(f_1) \odot T(f_2) \oplus_{\beta,S} T(f_1 \odot f_2) \odot \log (-\lambda) = T(T(f_1) \odot f_2) \oplus_{\beta,S} 
  T(f_1 \odot T(f_2)). $$
\endproof

We also check that linearity (in the ordinary sense) for the operator $\cT$ corresponds to
$\oplus_{\beta,S}$-linearity for $T$. For a semiring $\bS=C(X,\bT)$, we
extend $\cT$ to $\cR=C(X,\R)$ by requiring that $\cT(-e^{-\beta f}):=
- \cT(e^{-\beta f})$.

\begin{prop}\label{Tlinearbeta}
Let $\cT(e^{-\beta f}):= e^{-\beta T(f)}$, as in Theorem \ref{RBexp}. Then the operator $\cT$
is $\R$-linear if and only if the operator $T$ is $\oplus_{\beta,S}$-linear. 
\end{prop}

\proof We have
$$ \cT(e^{-\beta f_1} + e^{-\beta f_2})= \cT( e^{-\beta (-\beta^{-1}\log(e^{-\beta f_1}+e^{-\beta f_2}))}) $$
$$ = \cT(e^{-\beta( f_1\oplus_{\beta,S} f_2) }) =e^{-\beta (T(f_1\oplus_{\beta,S} f_2) )}. $$
We also have
$$ \cT(e^{-\beta f_1})+\cT(e^{-\beta f_2}) = e^{-\beta T(f_1)} + e^{-\beta T(f_2)} $$
$$ = e^{-\beta (-\beta^{-1} \log(e^{-\beta T(f_1)} + e^{-\beta T(f_2)}))} =
e^{-\beta T(f_1)\oplus_{\beta,S} T(f_2)}, $$
hence $\cT(e^{-\beta f_1} + e^{-\beta f_2})= \cT(e^{-\beta f_1})+\cT(e^{-\beta f_2}) $ if and
only if $T(f_1\oplus_{\beta,S} f_2)=T(f_1)\oplus_{\beta,S} T(f_2)$.
Moreover, for $\alpha \in \R^*_+$, we have 
$$ \cT( \alpha e^{-\beta f})= \cT( e^{-\beta (f-\beta^{-1} \log\alpha)}) =
e^{-\beta T(f-\beta^{-1} \log\alpha)}. $$
This agrees with
$$ \alpha \cT(e^{-\beta f}) = \alpha e^{-\beta T(f)} = e^{-\beta (T(f)-\beta^{-1} \log \alpha)} $$
if and only if, for all $f\in C(X,\R)$ and all $\alpha\in \R^*_+$ we have $T(f-\beta^{-1} \log\alpha)=
T(f)-\beta^{-1} \log \alpha$.
The two properties $T(f_1\oplus_{\beta,S} f_2)=T(f_1)\oplus_{\beta,S} T(f_2)$ and
$T(f+\lambda)= T(f)+\lambda$, for all $f, f_1, f_2\in C(X,\R)$ and all $\lambda\in \R$, are
equivalent to $\oplus_{\beta,S}$-linearity \eqref{oplusbetaSlinear}.
\endproof

\medskip
\subsection{Birkhoff factorization in thermodynamic semirings}

Let $\cH$ be a graded connected commutative Hopf algebra and
$\psi: \cH \to \bS_{\beta,S}$ satisfying $\psi(xy)=\psi(x)+\psi(y)$.

\begin{defn}\label{BPHZbeta}
Let $T: \bS_{\beta,S}\to \bS_{\beta,S}$ be a Rota--Baxter operator of weight $\lambda=+1$,
as in Definition \ref{RBplusbetaS}.
The Bogolyubov--Parashchuk preparation of $\psi$ is defined as
\begin{equation}\label{preparebeta}
\tilde\psi_{\beta,S}(x) =\psi(x) \oplus_{\beta,S} \bigoplus_{\beta,S} \psi_-(x') + \psi(x'') =
-\beta^{-1} \log \left(e^{-\beta \psi(x)} + \sum e^{-\beta (\psi_-(x') + \psi(x'')) } \right),
\end{equation}
where $\Delta(x)=x\otimes 1+1\otimes x + \sum x'\otimes x''$, and where
$\psi_-(x)=T\tilde\psi(x)$.
\end{defn}

\smallskip

\begin{rem}\label{limBP}{\rm
When $\beta \to \infty$, the Bogolyubov--Parashchuk preparation $\tilde\psi_{\beta,S}(x)$
converges to the preparation \eqref{BPpsi}. }
\end{rem}

\smallskip

\begin{lem}\label{BPHZphipsi}
Given $\psi: \cH \to \bS_{\beta,S}$ satisfying $\psi(xy)=\psi(x)+\psi(y)$, for all $x,y\in \cH$, let
$\phi_\beta (x):= e^{-\beta \psi(x)}$. Then $\phi_\beta(xy)=\phi_\beta(x)\phi_\beta(y)$, for all $x,y\in \cH$. The
Bogolyubov--Parashchuk preparation of $\psi$ satisfies $\tilde\phi_\beta(x)=e^{-\beta \tilde\psi(x)}$,
for all $x\in \cH$, where 
$$ \tilde\phi_\beta(x): = \phi_\beta(x) + \sum \cT(\tilde\phi_\beta(x')) \phi_\beta(x''), $$
with $\Delta(x)=x\otimes 1+1\otimes x + \sum x'\otimes x''$, and where the operator $\cT$
is defined by $\cT(e^{-\beta f}):= e^{-\beta T(f)}$ and $\cT(-e^{-\beta f}):=
- \cT(e^{-\beta f})$.
\end{lem}

\proof The multiplicativity of $\phi_\beta$ is evident. For the 
Bogolyubov--Parashchuk preparation, we inductively assume that for the
lower degree terms $\tilde\phi_\beta(x')= e^{-\beta \tilde\psi_\beta(x')}$.
Then the result follows from the relation between the operators $T$ and $\cT$.
\endproof

\smallskip

\begin{defn}\label{defstarbeta}
Given $\psi_1$ and $\psi_2$ from $\cH$ to $\bS_{\beta,S}$ with
$\psi_i(xy)=\psi_i(x)+\psi_i(y)$, we set
\begin{equation}\label{starbeta}
(\psi_1\star_\beta \psi_2)(x)=\bigoplus_{\beta,S} (\psi_1(x^{(1)}) + \psi_2(x^{(2)})), 
\end{equation}
where $\Delta(x)=\sum x^{(1)}\otimes x^{(2)}$.
\end{defn}

\smallskip

\begin{lem}\label{starbetastar}
Let $\phi_{i,\beta} (x)=e^{-\beta \psi_i(x)}$. Then 
$(\phi_{1,\beta}\star \phi_{2,\beta}) (x) = e^{-\beta (\psi_1\star_\beta \psi_2)(x)}$.
\end{lem}

\proof The usual product of Hopf algebra characters is given by
$$ (\phi_{1,\beta}\star \phi_{2,\beta}) (x) = \sum \phi_{1,\beta}(x^{(1)}) \, \phi_{2,\beta} (x^{(2)}) , $$
where $\Delta(x)=\sum x^{(1)}\otimes x^{(2)}$. This can be written equivalently as
$$ \sum e^{-\beta  ( \psi_{1,\beta}(x^{(1)}) + \phi_{2,\beta}) (x^{(2)})) } =
e^{-\beta ( \beta^{-1} \log ( \sum e^{-\beta  ( \psi_{1,\beta}(x^{(1)}) + \phi_{2,\beta} (x^{(2)})) } )} $$
$$ = e^{-\beta \bigoplus_{\beta,S} (\psi_1(x^{(1)}) + \psi_2(x^{(2)}))} = e^{-\beta (\psi_1\star_\beta \psi_2)(x)} . $$
\endproof

\smallskip

The construction of Birkhoff factorizations in thermodynamic semirings is then given by the following.

\begin{thm}\label{Birk1betaProp}
Let $T: \bS_{\beta,S} \to \bS_{\beta,S}$ be a $\oplus_{\beta,S}$-additive Rota--Baxter operator of 
weight $\lambda=+1$, in the sense of Definition \ref{RBplusbetaS}. 
Then there is a factorization $\psi_{\beta,+} =\psi_{\beta,-} \star_\beta \psi$, where $\psi_{\beta,\pm}$ are defined as
\begin{equation}\label{psibetamin}
\psi_{\beta,-}(x) = T (\tilde\psi_\beta(x)) = -\beta^{-1} \log \left(e^{-\beta T(\psi(x))} + \sum e^{-\beta T(\psi_-(x') + \psi(x'')) } \right)
\end{equation}
\begin{equation}\label{psibetaplus}
\psi_{\beta,+}(x) = -\beta^{-1} \log \left(e^{-\beta \psi_{\beta,-}(x)} + e^{-\beta \tilde\psi_\beta(x)} \right).
\end{equation}
The positive and negative parts of the Birkhoff factorization satisfy $\psi_{\beta,\pm}(xy)=\psi_{\beta,\pm}(x)+\psi_{\beta,\pm}(y)$.
\end{thm}

\proof By Lemma \ref{BPHZphipsi}, the statement is analogous to showing the existence of a factorization
$\phi_{\beta,+}=\phi_{\beta,-} \star_\beta \phi$ for $\phi_\beta (x)=e^{-\beta \psi(x)}$, with 
$\phi_{\beta,-}(x)=e^{-\beta \psi_{\beta,-}(x)}$ and $\phi_{\beta,+}(x)= e^{-\beta \psi_{\beta,+}(x)}$, and
such that $\phi_{\beta,\pm}(xy)=\phi_{\beta,\pm}(x) \phi_{\beta,\pm}(y)$. Such a factorization can be constructed
inductively by setting 
$$ \phi_{\beta,-}(x) = \phi_\beta(x) + \sum \phi_{\beta,-}(x') \phi_\beta(x''), $$
$$ \phi_{\beta,+}(x) = \phi_{\beta,-}(x) + \tilde\phi_\beta(x), $$
where, according to Lemma \ref{BPHZphipsi} and Propositions \ref{RBexp} and \ref{Tlinearbeta}, 
$$ \phi_{\beta,-}(x)=\cT (\tilde \phi_\beta(x))=e^{-\beta T(\tilde\psi_\beta(x))} =e^{-\beta \psi_{\beta,-}(x)}. $$
The multiplicative property for $\phi_{\beta,+}$ follows from that of $\phi_{\beta,-}$ and of $\phi_\beta$. Thus,
it suffices to show $\phi_{\beta,-}(xy)=\phi_{\beta,-}(x)\phi_{\beta,-}(y)$. We proceed as in the case of the
usual Birkhoff factorization, and identify the terms in 
$$ \phi_{\beta,-}(xy) = \cT (\phi_\beta(xy) + \sum \phi_{\beta,-}((xy)')\phi_\beta((xy)'') ) $$
with
$$ \cT ( \tilde\phi_\beta(x) \tilde \phi_\beta(y))  +\cT( \cT(\tilde\phi_\beta(x)) \tilde \phi_\beta(y) )+
\cT(\tilde\phi_\beta(x) \cT(\tilde \phi_\beta(y) )). $$
Using Proposition \ref{RBexp} and the resulting Rota--Baxter identity of weight $\lambda=+1$ for $\cT$, 
we identify this with 
$$ \cT( \tilde\phi_\beta(x)) \cT(\tilde \phi_\beta(y)) = \phi_{\beta,-}(x) \phi_{\beta,-}(y). $$
\endproof

\smallskip

\begin{rem}\label{limiBirkbeta0}{\rm
In the limit when $\beta \to \infty$, the Birkhoff factorization of Theorem \ref{Birk1betaProp} converges to the
Birkhoff factorization of Theorem \ref{algrensemi}.
}\end{rem}

\bigskip
\section{von Neumann entropy and Rota--Baxter structures}\label{vonNeumannSec}

We now consider again the case of matrices. Recall from Theorem 1.2.8 of \cite{Guo} that if $\cR$ is a commutative
$\R$-algebra, endowed with a Rota--Baxter operator $\cT$ of weight $\lambda$, then $\cT$ induces a Rota--Baxter
operator (which we still denote $\cT$), of the same weight, on the ring of matrices $M_n(\cR)$, by applying $\cT$
coordinate-wise, $\cT(A)=( \cT(a_{ij}) )$, for $A=(a_{ij})$.

\smallskip

\begin{prop}\label{dirsumRB}
Let $\cR$ be a commutative $\R$-algebra and $\bS$ be a min-plus semiring, related by the property that,
for $A\in M_n(\bS)$, the matrix $e^{-\beta A} \in M_n(\cR)$. Let $\cD\subset M_n(\cR)$ denote 
the set of matrices of the form $e^{-\beta A}$, for $A\in M_n(\bS)$. Let $\cT$ be a Rota--Baxter operator 
of weight $+1$ on $\cR$, and let $(M_n(\cR), \cT)$ be Rota--Baxter structure described above.
Then setting $\cT(e^{-\beta A})=e^{-\beta T(A)}$ defines an operator $T: M_n(\bS)\to M_n(\bS)$
that satisfies the following type of Rota--Baxter identity
\begin{equation}\label{RBmatrix}
\Tr^\oplus_{\beta,\cN}(T(A)\boxplus T(B)) = \Tr^\oplus_{\beta,\cN}(T(T(A)\boxplus B)) \oplus_{\beta,S}
\Tr^\oplus_{\beta,\cN}(T(A\boxplus T(B))) \oplus_{\beta,S} \Tr^\oplus_{\beta,\cN}(T(A)\boxplus T(B)),
\end{equation}
where $\cN$ is the von Neumann entropy, $S$ is the Shannon entropy, and $\boxplus$ denotes 
the direct sum of matrices.
\end{prop}

\proof
The Rota--Baxter identity for $\cT$ on $\cR$ gives
\begin{eqnarray*}
 \cT(\Tr(e^{-\beta A})) \cT(\Tr(e^{-\beta B})) & = & \cT(\cT(\Tr(e^{-\beta A}))  \Tr(e^{-\beta B}))
+\cT( \Tr(e^{-\beta A}) \cT(\Tr(e^{-\beta B})) \\[2mm] & + & \cT( \Tr(e^{-\beta A}) \Tr(e^{-\beta B})). 
\end{eqnarray*}
Notice that the induced Rota--Baxter structure on $M_n(\cR)$ satisfies $\cT(\Tr(A))=\Tr(\cT(A))$.
Thus, using this fact together with $\cT(e^{-\beta A})=e^{-\beta T(A)}$, we can rewrite the above as
\begin{eqnarray*}
\Tr( e^{-\beta T(A)}) \Tr(e^{-\beta T(B)}) & = & \cT( \Tr(e^{-\beta T(A)}) \Tr(e^{-\beta B}))
+ \cT( \Tr(e^{-\beta A}) \Tr(e^{-\beta T(B)}) ) \\[2mm] & + & \cT( \Tr(e^{-\beta A}) \Tr(e^{-\beta B})). 
\end{eqnarray*}
We can then identify the products of traces with the trace of the tensor product of matrices,
which gives
$$ \Tr( e^{-\beta T(A)} \otimes e^{-\beta T(B)}) = \cT( \Tr(e^{-\beta T(A)} \otimes e^{-\beta B}))
+ \cT( \Tr(e^{-\beta A} \otimes e^{-\beta T(B)}) ) + \cT( \Tr(e^{-\beta A} \otimes e^{-\beta B})). $$
Moreover, for matrix exponentials, $\exp(A)\otimes\exp(B)=\exp(A \boxplus B)$, where here $\boxplus$
is the direct sum of matrices. Thus, we obtain 
\begin{eqnarray*}
\Tr ( e^{-\beta(T(A)\boxplus T(B))}) & = & \Tr(e^{-\beta(T(T(A)\boxplus B))}) +
\Tr(e^{-\beta(T(A\boxplus T(B))}) \\[2mm] & + & \Tr(e^{-\beta(T(A)\boxplus T(B))}). 
\end{eqnarray*}
This then gives
\begin{eqnarray*}
-\beta^{-1} \log\Tr ( e^{-\beta(T(A)\boxplus T(B))}) & = & -\beta^{-1} \log ( \Tr(e^{-\beta(T(T(A)\boxplus B))}) +
\Tr(e^{-\beta(T(A\boxplus T(B))}) \\[2mm] & + & \Tr(e^{-\beta(T(A)\boxplus T(B))}) ),
\end{eqnarray*}
or equivalently
$$ 
\Tr^\oplus_{\beta,\cN}(T(A)\boxplus T(B))  =  -\beta^{-1} \log \left( e^{-\beta \Tr^\oplus_{\beta,\cN}(T(T(A)\boxplus B))}
+ e^{-\beta  \Tr^\oplus_{\beta,\cN}( T(A\boxplus T(B)) ) } 
 +  e^{-\beta  \Tr^\oplus_{\beta,\cN}(T(A)\boxplus T(B)) } \right), 
$$
hence \eqref{RBmatrix} follows.
\endproof

\bigskip
\section{Rota--Baxter structures of weight one: thermodynamics and Witt rings}\label{WittSec}

We analyze here some examples for Rota--Baxter operators of weight $+1$ on thermodynamic
semirings, derived from classical examples of weight-one Rota--Baxter algebras. We also show
that the same examples of weight-one Rota--Baxter algebras can be used to induce Rota--Baxter
structures on Witt rings. We interpret the effect of the resulting Rota--Baxter operators applied to
zeta functions of varieties, regarded as elements of Witt rings, as in \cite{Rama}.

\medskip
\subsection{Partial sums}\label{parsumSec}

Consider the $\R$-algebra $\cR$ of $\R$-valued sequences  
$a=(a_1,a_2,a_3,\cdots)=(a_n)_{n=1}^\infty$, with coordinate-wise addition and multiplication,
and let $\cT : \cR \to \cR$ be the linear operator that maps the sequence $(a_1,a_2, a_3, \cdots, a_n, \cdots)$
to $(0,a_1,a_1+a_2, \cdots, \sum_{k=1}^{n-1} a_k, \cdots)$. The operator $\cT$ is a Rota--Baxter operator
of weight $+1$, see Example 1.1.6 of \cite{Guo}.

\smallskip

\begin{lem}\label{Tseq}
The Rota--Baxter algebra $(\cR,\cT)$ of weight $+1$ described above determines a Rota--Baxter structure
of weight $+1$ on the thermodynamic semi-rings $\bS_{\beta,S}$ of functions 
$f: \N \to \bT=\R\cup \{\infty \}$, with the pointwise operations $\oplus_{\beta,S}$ and $\odot$, with
Rota--Baxter operator
\begin{equation}\label{RBseqbeta}
(T f)(n) = {\bigoplus_{\beta,S}}_{k=1,\ldots,n-1}  f(k) ,
\end{equation}
for $n\geq 2$ and $(T f)(1)=\infty$.
\end{lem}

\proof
For $\cR$ as above, let $\cD\subset \cR$ be the subset of sequences with values in $\R_+$,
which we can write as $a_n=e^{-\beta c_n}$, when $a_n>0$ and zero otherwise. 
We have $(\cT a)_1=0$ and $(\cT a)_n=\sum_{k=1}^{n-1} a_k$ for $n\geq 2$.
Define $(T c)_n =- \beta^{-1} \log (\cT a)_n$, so that 
$$ (T c)_n = \left\{ \begin{array}{ll} \infty & n=1 \\[2mm]
-\beta^{-1} \log \left( \sum_{k=1}^{n-1} e^{-\beta c_k} \right) & n \geq 2 .
\end{array} \right. $$
\endproof

\medskip
\subsection{$q$-integral}

Let $\cR=\R[[t]]$ be the ring of formal power series with real coefficients. Let $\cT$ be the linear operator 
$(\cT \alpha)(t)=\sum_{k=1}^\infty \alpha(q^n t)$, for $q$ not a root of unity. The operator $\cT$ is a
Rota--Baxter operator of weight $+1$. The operator $\cT$ maps a single power $t^n$ to $q^n t^n/(1-q^n)$,
hence it restricts to a Rota--Baxter operator of weight $+1$ on the subring of polynomials $\R[t]$,
see Example 1.1.8 of \cite{Guo}.

\smallskip

\begin{lem}
Let $\bS$ be the thermodynamic semiring of formal power series $\bS_{\beta,S}=\R[[t]] \cup \{ \infty \}$
with the operations $(\gamma_1\oplus_{\beta,S} \gamma_2)(t)
=-\beta^{-1}\log( e^{-\beta \gamma_1(t)} + e^{-\beta \gamma_2(t)})$ and with $(\gamma_1\odot \gamma_2)(t)
=\gamma_1(t)+\gamma_2(t)$. Then the Rota--Baxter algebra of weight $+1$, given by the data $(\cR,\cT)$
described above, induces a Rota--Baxter structure of weight $+1$ on $\bS_{\beta,S}$ by
\begin{equation}\label{RBqseriesbeta}
(T \gamma)(t) = {\bigoplus_{\beta,S}}_{k=1}^\infty  \gamma(q^k t).
\end{equation}
\end{lem}

\proof Let $\cD\subset \cR$ be the subset of formal series with $a_0=1$, that is,
$\cD= 1+t \R[[t]]$. Then for $\alpha \in \cD$ and $\gamma(t)=\log\alpha(t)$,we define
an operator $T$ by the relation $\cT(e^{-\beta \gamma(t)})=e^{-\beta (T\gamma)(t)}$.
This gives $e^{-\beta (T \gamma(t))}=\sum_{k=1}^\infty e^{-\beta \gamma(q^k t)}$,
that is,
$$ (T\gamma)(t)= -\beta^{-1} \log\left(  \sum_{k=1}^\infty e^{-\beta \gamma(q^k t)} \right) = 
{\bigoplus_{\beta,S}}_{k=1}^\infty  \gamma(q^k t). $$
\endproof

\medskip
\subsection{Rota--Baxter structures on Witt rings}

For a commutative ring $R$, the Witt ring $W(R)$ can be identified with the
set of formal power series with $a_0=1$, that is, the set $1+ t R[[t]]$ with the Witt 
addition given by the usual product of formal power series and the Witt multiplication 
$\star$ uniquely determined by the rule
$$ (1-at)^{-1} \star (1-bt)^{-1} = (1-ab t)^{-1}, $$
for $a,b \in R$. There is an injective ring homomorphism $g: W(R)\to R^\N$,
$g(\alpha)=(\alpha_1,\alpha_2,\ldots,\alpha_r,\ldots)$, where the addition and
multiplication operations on $R^\N$ are component-wise. The sequence 
$g_n(\alpha)=\alpha_n$ is known as the ``ghost coordinates" of $\alpha$.
Upon writing elements of the Witt ring $W(R)$ in the exponential form
$$ \exp\left(\sum_{r\geq 1} \alpha_r \frac{t^r}{r} \right), $$
one sees that the the ghost coordinates are the coefficients of
$$ t \frac{1}{\alpha} \frac{d\alpha}{dt} = \sum_{r\geq 1} \alpha_r t^r. $$

\smallskip

\begin{lem}\label{RBghost}
A linear operator $\cT : R^\N \to R^\N$ is a Rota--Baxter operator of weight $\lambda$ 
on $R^\N$ if and only if the operator $\cT_W$ defined on the Witt ring $W(R)$ so that,
when taking ghost components $g(\cT_W (\alpha))=\cT(g(\alpha))$ is a
Rota--Baxter operator of weight $\lambda$ on $W(R)$.
\end{lem}

\proof  the Rota--Baxter identity for $\cT_W$ is of the form
$$ \cT_W(\alpha_1)\star \cT_W(\alpha_2)  =  \cT_W( \alpha_1 \star \cT_W(\alpha_2)) +_W
\cT_W (\cT_W(\alpha_1)\star \alpha_2) +_W \lambda \star \cT_W (\alpha_1 \star \alpha_2)) $$
with $+_W$ the sum in $W(R)$. When taking ghost components, this gives
$$ g( \cT_W(\alpha_1)\star \cT_W(\alpha_2)) ) = g( \cT_W( \alpha_1 \star \cT_W(\alpha_2))) +
g(\cT_W (\cT_W(\alpha_1)\star \alpha_2)) +\lambda \, g (\cT_W (\alpha_1 \star \alpha_2)) $$
which gives the Rota--Baxter identity for $\cT$, 
$$ \cT(g(\alpha_1)) \cT(g(\alpha_2))=  \cT ( g( \alpha_1) \cT(g(\alpha_2)))  +
\cT( \cT(g(\alpha_1)) g(\alpha_2)) +\lambda \cT( g (\alpha_1) g(\alpha_2)). $$
The injectivity of the ghost map shows we can run the implication backward.
\endproof

\smallskip

In addition to the Witt product $\star$ of the Witt ring $W(R)$, which corresponds to
the coordinate-wise product of the ghost components, one can introduce a convolution
product on $W(R)$, which is induced by the power-series product of the ghost maps.

\begin{defn}\label{Wittconvol}
For $\alpha,\gamma \in W(R)$, with
$\alpha=\exp(\sum_{r\geq 1} \alpha_r t^r/r)$ and $\gamma=\exp(\sum_{r\geq 1} \gamma_r t^r/r)$,
the convolution product is given as
\begin{equation}\label{dastdef}
 \alpha \circledast \gamma := \exp\left( \sum_{n\geq 1} (\sum_{r+\ell=n} \alpha_r \gamma_\ell)  \frac{t^n}{n} \right). 
\end{equation}
\end{defn}

Notice that $\alpha \circledast \gamma$ is defined so that the ghost
$g(\alpha \circledast \gamma)=\sum_{n\geq 1} \sum_{r+\ell=n} \alpha_r \gamma_\ell \, t^n$ is the
product as power series $g(\alpha) \bullet g(\gamma)$ of the 
ghosts $g(\alpha)=\sum_{r\geq 1} \alpha_r t^r$ and $g(\gamma)=\sum_{r\geq 1} \gamma_r t^t$.

\smallskip

\begin{lem}\label{convolRBop}
A linear operator $\cT : R[[t]] \to R[[t]]$ is a Rota--Baxter operator of weight $\lambda$ if and
only if the operator $\cT_W: W(R) \to W(R)$ defined by $g(\cT_W(\alpha))=\cT(g(\alpha))$
satisfies the Rota--Baxter identity of weight $\lambda$ with respect to the convolution product \eqref{dastdef},
\begin{equation}\label{convolRBid}
\cT_W(\alpha_1) \circledast \cT_W(\alpha_2) = \cT_W (\alpha_1 \circledast \cT_W(\alpha_2)) +_W \cT_W (\cT_W(\alpha_1)
\circledast \alpha_2) +_W \lambda \cT_W(\alpha_1 \circledast \alpha_2),
\end{equation}
where $+_W$ is the addition in $W(R)$.
\end{lem}

\proof After composing with the ghost map, \eqref{convolRBid} gives
$$ \cT(g(\alpha_1)) \bullet \cT(g(\alpha_2)) = \cT(g(\alpha_1)\bullet \cT(g(\alpha_2)))+
\cT(\cT(g(\alpha_1))\bullet g(\alpha_2)) + \lambda \, \cT(g(\alpha_1)\bullet g(\alpha_2)), $$
where $\bullet$ denotes the product as formal power series. This is the
Rota--Baxter identity for $\cT$ on $R[[t]]$. The injectivity of the ghost map shows
the two conditions are equivalent.
\endproof

\medskip

We consider then the example of Rota--Baxter operator of weight one given by partial sums.

\begin{prop}\label{RBWitt1}
Let $\cR=R^\N$ with the Rota--Baxter operator of weight $+1$ given by 
$$ \cT: (a_1,a_2,\ldots, a_n, \ldots) \mapsto (0, a_1, a_1+a_2, \ldots, \sum_{k=1}^{n-1} a_k, \ldots). $$
The resulting Rota--Baxter operator $\cT_W$  of weight $+1$ on the Witt ring $W(R)$ is given by 
convolution product with the multiplicative unit (for the usual Witt product) ${\mathbb I}=(1-t)^{-1}$ of $W(R)$,
\begin{equation}\label{TWconvol}
\cT_W (\alpha) = \alpha \circledast {\mathbb I}.
\end{equation}
\end{prop}

\proof According to Lemma \ref{RBghost}, the Rota--Baxter operator $\cT_W$ on $W(R)$ is given by
$$ \cT_W(\alpha) = \exp\left( \sum_{n\geq 2} \sum_{k=1}^{n-1} \alpha_k \, \frac{t^n}{n} \right), $$
for $\alpha=\exp(\sum_{n\geq 1} \alpha_n t^n/n )$. The ghost of $\alpha$ is given by
$g(\alpha)=\sum_{n\geq 1}\alpha_n t^n$ and we can identify the series
$$ \sum_{n\geq 2} \sum_{k=1}^{n-1} \alpha_k \, t^n = g(\alpha) \bullet \frac{t}{1-t}, $$
where $\bullet$ denotes the product of formal power series. The identity
$$ - t \frac{d}{dt} \log (1-t) = \frac{t}{1-t} $$
then shows that we can identify the above with the product of power series $g(\alpha)\bullet g({\mathbb I})$,
hence by construction $\cT_W(\alpha)=\alpha \circledast {\mathbb I}$.
\endproof

\medskip

We also consider the example of the weight-one Rota--Baxter operator on power series given by the $q$-integral.

\begin{prop}\label{RBWitt2}
Let $\cR=R[[t]]$ with the Rota--Baxter operator $\cT_q$ of weight $+1$ given by the $q$-integral (where $q\in R$ is not
a root of unity). Then the operator $\cT_{W,q}$ on $W(R)$ defined by $g(\cT_{W,q}(\alpha))=\cT_q(g(\alpha))$ is a
Rota--Baxter operator of weight one with respect to the convolution product \eqref{dastdef}. It is explicitly
given by $\cT_W(\alpha)(t)=\prod_{k\geq 1} \alpha(q^k t)$.
\end{prop}

\proof The operator $\cT_{W,q}$ acts as 
$$ \cT_{W,q}(\exp(\sum_{r\geq 1} \alpha_r \frac{t^r}{r}))=\exp(\sum_{r\geq 1} \sum_{k\geq 1} \alpha_r \frac{q^{kr} t^r}{r}))
= \prod_{k\geq 1} \exp(\sum_{r\geq 1} \alpha_r \frac{(q^k t)^r}{r}). $$
Notice that the product $\prod_k \alpha(q^k t)$, which is the product as power series, is the {\em addition}
in the Witt ring $W(R)$, so the operator $\cT_W$ has the same form as the $q$-integral operator $\cT$, simply
replacing the sum in $R[[t]]$ with the sum in $W(R)$. By Lemma \ref{convolRBop}, $\cT_{W,q}$ satisfies the identity 
$$ \cT_{W,q}(\alpha_1) \circledast \cT_{W,q}(\alpha_2) = \cT_{W,q} (\alpha_1 \circledast \cT_{W,q}(\alpha_2)) 
+_W \cT_{W,q} (\cT_{W,q}(\alpha_1)\circledast \alpha_2) +_W \cT_{W,q}(\alpha_1 \circledast \alpha_2). $$
\endproof

\medskip
\subsection{Applications to zeta functions}

For varieties (or schemes) over finite fields, the Hasse--Weil zeta function is given by
$$ Z(X,t)= \exp\left( \sum_{r\geq 1} \# X(\F_{q^r}) \, \frac{t^r}{r} \right). $$
Equivalently, it can be written as
$$ Z(X,t)=  \prod_{r\geq 1} (1-t^r)^{-a_r(X)} = \prod_{x\in X_{cl}} (1-t^{\deg(x)})^{-1}, $$
where $a_r(X)=\# \{ x\in X_{cl} \,|\, \, [k(x):\F_q]=r \}$ and $X_{cl}$ is the set of closed points of $X$.
The zeta function satisfies the properties
$$ Z(X\sqcup Y, t)=Z(X,t) Z(Y,t), $$
for a disjoint union $X\sqcup Y$ and
$$ Z(X \times Y, t) = Z(X,t) \star Z(Y,t), $$
where $\star$ is the product in the Witt ring. Thus, it is natural to consider zeta functions
of varieties as elements of a Witt ring, \cite{Rama}. In particular, this means that we can
apply the Rota--Baxter operators on Witt rings described above to zeta functions of varieties.

\smallskip

\begin{cor}\label{RBZeta1}
Let $\cT_W$ be the Rota--Baxter operator of Proposition \ref{RBWitt1}. For $X$ a variety (or scheme)
over $\F_q$, 
$$ \cT_W(Z(X,t))= Z(X,t) \circledast Z(\Spec(\F_q),t). $$
\end{cor}

\proof This is immediate from Proposition \ref{RBWitt1}, since $Z(\Spec(\F_q),t)=\exp(\sum_{r\geq 1} \frac{t^r}{r})=(1-t)^{-1}$.
\endproof

\smallskip

Recall that the Grothendieck ring of varieties (or schemes of finite type) over $\F_q$ is generated by isomorphism classes 
$[X]$ with the inclusion-exclusion relation $[X]=[Y]+[X\smallsetminus Y]$ for closed $Y\subset X$ and the product
$[X\times Y]=[X] \, [Y]$. The zeta function $Z(X,t)=Z([X],t)$ factors as a ring homomorphism from the Grothendieck ring
to the Witt ring. In the Grothendieck ring, the Lefschetz motive is the class of the affine line $\bL=[\A^1]$. In the
theory of motives it is customary to localize the Grothendieck ring by inverting the Lefschetz motive. The Tate
motive is the formal inverse $\bL^{-1}$. For the Rota--Baxter structure of Proposition \ref{RBWitt2} we then have the following.

\begin{cor}\label{RBZeta2}
Let $X$ be a variety (or scheme) over $k=\F_q$ and $[X]$ its Grothendieck class. Consider the Rota--Baxter structure of 
Proposition \ref{RBWitt2} with Rota--Baxter operator $\cT_{W,q}$ or $\cT_{W,q^{-1}}$. These
give
$$ \cT_{W,q}(Z(X,t))=\prod_{k\geq 1} Z([X] \,\bL^k,t), \ \ \ \  \cT_{W,q^{-1}}(Z(X,t))=\prod_{k\geq 1} Z([X] \,\bL^{-k},t), $$
where $\bL$ is the Lefschetz motive and $\bL^{-1}$ is the Tate motive.
\end{cor}

\proof In the case of the Lefschetz motive we have $Z(X, q^kt)=Z(X \times \A^k, t)=Z([X]\bL^k, t)$. 
In the case of the Tate motive, we do not have the geometric space replacing $X \times \A^k$, but
the property that the zeta function is a ring homomorphism from the Grothendieck ring to the Witt ring
gives $Z(X, q^{-k} t)=Z([X]\bL^{-k}, t)$.
\endproof

\smallskip

\begin{cor}\label{RBZeta3}
Let $\cT_{W,q}$ and $\cT_{W,q^{-1}}$ be as above. Then the operators $\tilde\cT_{W,q^{\pm 1}}:=
-_W id -_W \cT_{W,q^{\pm 1}}$ are also Rota--Baxter operators of weight $+1$. 
For $X$ a variety over $\F_q$, they give
$$ \tilde\cT_{W,q^{\pm 1}}(Z(X,t))=\prod_{k\geq 0} Z([X] \bL^{\pm k}, t)^{-1}. $$
\end{cor}

\proof It is a simple general fact that, of $\cT$ is a Rota--Baxter operator of weight $+1$ then
$\tilde\cT=-id -\cT$ is also a Rota--Baxter operator of weight $+1$.  Thus, the operators $\tilde\cT_{W,q^{\pm 1}}$
satisfy the identity \eqref{convolRBid} with $\lambda=+1$. The explicit expression for 
$\tilde\cT_{W,q^{\pm 1}}(Z(X,t))$ then follows exactly as in Corollary \ref{RBZeta2}.
\endproof

\bigskip
\section{Rota--Baxter operators of weight $-1$}\label{minoneSec}

We have seen in the previous sections how to construct Birkhoff factorizations in min-plus
semirings and their thermodynamical deformations, based on the use of Rota--Baxter operators
of weight $\lambda=+1$.

\smallskip

In the usual setting of Birkhoff factorizations in perturbative quantum field theory, \cite{CoKr}, \cite{EFGuo},
\cite{EFGK}, one constructs the Birkhoff factorization using a Rota--Baxter operator $\cT$ of weight $\lambda=-1$
by setting  $\phi_-(x) = -\cT (\tilde\phi(x))$, where $-\cT$ is a Rota--Baxter operator of weight $+1$.

\smallskip

In the semiring setting, one cannot proceed in the same way. However, it is still possible to construct 
Birkhoff factorizations from Rota--Baxter operator of weight $\lambda=-1$, under some additional
conditions on the operator.

\smallskip

Let $\bS=C(X,\bT)$ with the pointwise min-plus operations $\oplus$ and $\odot$.

\begin{prop}\label{algrensemimin1}
Let $\psi: \cH \to \bS$ be a min-plus character, and let $T: \bS \to \bS$ be
a Rota-Baxter operator of weight $-1$, in the sense of Definition \ref{RBsemimin}.
Then there is a Birkhoff factorization $\psi_+ = \psi_- \star \psi$.
If, moreover, the Rota-Baxter operator $T$ satisfies $T(f_1+f_2)\geq T(f_1)+T(f_2)$,
then $\psi_-$ and $\psi_+$ are also min-plus characters.
\end{prop}

\proof As in the case of Theorem \ref{algrensemi}, we define the two sides of the
factorization as $\psi_-(x) := T(\tilde\psi(x))$ and $\psi_+(x) :=(\psi_- \star \psi)(x)=\min\{ \psi_-(x), \tilde\psi(x) \}$,
where the preparation $\tilde\psi(x)$ is defined as in \eqref{BPpsi}. To show that $\psi_-(xy)=\psi_-(x)+\psi_-(y)$,
we again list the terms $(xy)'$ and $(xy)''$ as in Theorem \ref{algrensemi} and obtain
$$\psi_-(xy)=\min \{  T(\tilde\psi(x) + \tilde\psi(y)), T(T(\tilde\psi(x))+\tilde\psi(y)), T(\tilde\psi(x)+T(\tilde\psi(y))) \}. $$
The Rota--Baxter identity of weight $-1$ for the operator $T$ then gives
$$ \psi_-(xy)= \min\{ T(\tilde\psi(x) + \tilde\psi(y)), T(\tilde\psi(x)) + T(\tilde\psi(y))) \}. $$
If the operator $T$ satisfies $T(f_1+f_2)\geq T(f_1)+T(f_2)$, for all $a,b \in \bS$, we then
have
$$ \psi_-(xy)= T(\tilde\psi(x)) + T(\tilde\psi(y))) = \psi_-(x)+ \psi_-(y). $$
\endproof

\smallskip

The following observations show that there are choices of semiring
Rota-Baxter operators satisfying $T(a+b)\geq T(a)+T(b)$.

\smallskip

\begin{prop}\label{ProjLinRB}
For $X$ a compact Hausdorff space, 
let $\bS=C(X,\R)$, with the pointwise $\oplus, \odot$ operations. 
Let $T: \bS \to \bS$ be an idempotent linear operator (in the ordinary sense) 
on the underlying algebra $C(X,\R)$ with respect to the
usual additive structure on $C(X,\R)$. Then $T$ is (trivially) a semiring Rota-Baxter operator
satisfying the hypothesis of Theorem \ref{algrensemi}.
\end{prop}

\proof Since $T$ is linear, it satisfies $T(a+b)=T(a)+T(b)$ and the
Rota-Baxter relation simply becomes $T(a+b)=\min\{ T^2(a)+T(b), T(a)+T^2(b) \}$,
which is certainly satisfied if $T$ is idempotent, $T^2=T$, since the right-hand-side
is then also equal to $T(a)+T(b)$.
\endproof

\smallskip

\begin{ex}\label{exCantor} {\rm 
Let $\bS=C(X,\bT)$, as in Proposition \ref{ProjLinRB}, where $X$ is a totally disconnected
compact Hausdorff space (a Cantor set). Then for any clopen subset $Y\subset X$ the
operator $T=T_Y : \bS \to \bS$ given by ordinary multiplication by the characteristic function
of $Y$, $T_Y : f(x) \mapsto \chi_Y(x) f(x)$ is a semiring Rota--Baxter operator satisfying the
conditions of Theorem \ref{algrensemi}.}
\end{ex}

\smallskip

\begin{rem}\label{RBnonlin} {\rm
Examples of semiring Rota--Baxter operators satisfying the conditions of Theorem \ref{algrensemi},
but not arising from linear operators in the usual sense, can be constructed using idempotent
superadditive operators. These occur, for instance, in potential theory: we refer the reader to
\S 10 of \cite{ArsLeu} for some relevant constructions.}
\end{rem}

\bigskip
\subsection{Other forms of Birkhoff factorization of weight $-1$ in min-plus semirings}

We consider here a different possible way of obtaining Birkhoff factorization for min-plus semirings,
using Rota--Baxter operators of weight $-1$. This method involves
the use of two related Rota--Baxter operators, generalizing the
roles of the operators $\cT$ and $1-\cT$ in the original commutative algebra case.
As in the case of Proposition \ref{algrensemimin1}, we need superadditivity conditions
on these operators to obtain that the parts of the factorization are still min-plus character.

\smallskip

\begin{defn}\label{Birk2minplus}
Let $(\bS,\oplus,\odot)$ be a min-plus semiring. Let $T: \bS \to \bS$ and 
$\tilde T: \bS \to \bS$ be $\oplus$-additive Rota--Baxter operators  of weight $-1$
(as in Definition \ref{RBsemimin}) satisfying the relations
\begin{equation}\label{minTTtil}
T \alpha  =\alpha \oplus \tilde T \alpha, \ \ \  \forall \alpha \in \bS
\end{equation} 
\begin{equation}\label{TtildeTrel}
\tilde T(\alpha \odot \beta) \oplus  
\tilde T(\alpha) \odot \tilde T(\beta) =  \tilde T( T(\alpha) \odot \beta \oplus \alpha \odot T(\beta)), 
\ \ \  \forall \alpha,\beta\in \bS.
\end{equation}
Given a min-plus character 
$\psi: \cH \to \bS$,  a $(T,\tilde T)$-Birkhoff factorization of $\psi$ is given by the pair
\begin{equation}\label{psimRBstrong}
\psi_-(x) = T\tilde\psi(x)= T(\psi(x) \oplus \bigoplus_{(x',x'')} \psi_-(x') \odot \psi(x'') )
\end{equation}
\begin{equation}\label{psipRBstrong}
\psi_+(x)=\tilde T \tilde\psi(x)=\tilde T(\psi(x) \oplus \bigoplus_{(x',x'')} \psi_-(x') \odot \psi(x'') ),
\end{equation}
where the $\oplus$-sums are over pairs $(x',x'')$ in the non-primitive part of the
coproduct of $\cH$, $\Delta(x)=x\otimes 1+1\otimes x + \sum x' \otimes x''$.
\end{defn}

\smallskip

\begin{prop}\label{stBirkprop}
Let the data $(\bS,\oplus,\odot)$, with $T: \bS \to \bS$ and 
$\tilde T: \bS \to \bS$, be as in Definition \ref{Birk2minplus}. 
If the operators $T$ and $\tilde T$ are superadditive then 
the resulting pieces $\psi_\pm$ of the factorization are min-plus characters. 
Moreover,  the terms $\psi_\pm$ of the $(T,\tilde T)$-Birkhoff factorization
are related by 
\begin{equation}\label{psipmmin}
\psi_- (x)=\min\{ \tilde\psi(x), \psi_+(x) \}=\min\{ \psi_+ \star \psi(x), \tilde\psi(x')+\psi(x'') \}. 
\end{equation}
\end{prop}

\proof We first need to check that $\psi_\pm (xy)=\psi_\pm(x) +\psi_\pm(y)$. 
The case of $\psi_-$ is proved as in Theorem \ref{algrensemi}.
In the case of $\psi_+$, by proceeding as in Theorem \ref{algrensemi}, we see that
$$ \psi_+(xy)= \tilde T \min \{ T(\tilde\psi(x)) + \tilde\psi(y), \tilde\psi(x)+T(\tilde \psi(y)),
\tilde\psi(x)+\tilde\psi(y) \}. $$
Using the $\oplus$-additivity (monotonicity) of $\tilde T$ and \eqref{TtildeTrel} we
write the above as
$$ \psi_+(xy)=\min \{ \tilde T(\tilde\psi(x)) + \tilde T(\tilde\psi(y)), \tilde T(\tilde\psi(x)+\tilde\psi(y)) \}. $$
If $\tilde T$ is subadditive, the minimum is equal to 
$$ \psi_+(xy)= \tilde T(\tilde\psi(x)) + \tilde T(\tilde\psi(y))= \psi_+(x) + \psi_+(y). $$
Thus, both sides of the factorization satisfy $\psi_\pm (xy)=\psi_\pm(x) +\psi_\pm(y)$. 
We have
$$ \psi_-(x) = T(\tilde\psi(x))= \min \{ T(\psi(x)), T(T(\tilde\psi(x'))+\psi(x'')) \}. $$
The identity \eqref{minTTtil} implies that we have 
$$ \psi_-(x) = \min\{ \tilde\psi(x), \psi_+(x) \} $$
which we write also as $\min\{ \psi(x), \psi_+(x), \psi_-(x')+\psi(x'') \}$.
We then use \eqref{minTTtil} and rewrite $\psi_-(x')=\min \{ \tilde\psi(x'), \tilde T (\tilde\psi(x')) \}$.
Thus we obtain $\psi_-(x) = \min\{ \psi(x), \psi_+(x), \psi_+(x')+\psi(x''), \tilde\psi(x')+\psi(x'')\}$,
where $\min\{ \psi(x), \psi_+(x), \psi_+(x')+\psi(x'') \}$ is the convolution product
$(\psi_+\star \psi)(x)$, hence the statement follows.
\endproof

\bigskip

\section{Min-plus characters and thermodynamics}\label{examplesSec}

We consider here some examples of min-plus characters $\psi:\cH \to \bS$, satisfying
$\psi(xy)=\psi(x)+\psi(y)$. We focus on the case where $\cH$ a Hopf algebra of graphs,
namely the commutative algebra generated by connected finite graphs with coproduct
$$ \Delta(\Gamma)=\Gamma\otimes 1 + 1\otimes \Gamma + \sum_{\gamma \subset \Gamma} \gamma \otimes \Gamma/\gamma. $$

\medskip
\subsection{Inclusion--exclusion functions on graphs}

We consider real valued functions $\tau$ on a set of graphs, that satisfy an
inclusion-exclusion property. Namely, if $\Gamma =\Gamma_1 \cup \Gamma_2$
with intersection $\gamma =\Gamma_1 \cap \Gamma_2$, then
\begin{equation}\label{tauGamma}
\tau(\Gamma) = \tau(\Gamma_1) + \tau(\Gamma_2) - \tau(\gamma).
\end{equation}
Examples of such functions can be constructed by assigning
a ``cost function" to the sets of vertices and edges 
of a graph. Let $F_E=\{ f_e\,:\, e\in E(\Gamma) \}$ and $F_V=\{ f_v \,:\, v\in V(\Gamma) \}$.
Then setting $\tau(\Gamma)= \sum_{v\in V(\Gamma)} f_v + \sum_{e\in E(\Gamma)} f_e$
gives a function that satisfies inclusion-exclusion. (One of the sums may be trivial if one 
only assigns vertex or edge labels.) 

\smallskip

In particular, for a disjoint union $\Gamma = \Gamma_1 \sqcup \Gamma_2$ we 
have $\tau(\Gamma)=\tau(\Gamma_1)+\tau(\Gamma_2)$, hence we can view
such a function $\tau$ as a morphism $\tau: \cH \to \bT$, where $\cH$ is the
Hopf algebra of graphs and $\bT$ is the tropical semiring, satisfying
$\tau(xy)=\tau(x)+\tau(y)$, hence it defines a min-plus character.
The function $\tau$ obtained as above may depend on a set of parameters, so that
the $f_e$ and $f_v$ are functions of these parameters, so that we can think
of $\tau:\cH \to \bS$ as a min-plus character to some min-plus semiring of functions.

\medskip
\subsection{Examples from computation}

Following \S 4.6 of \cite{Man1}, we consider a Hopf algebra of ``flow charts" for
computation, namely graphs endowed with acyclic orientations, so that the flow through the graph,
from the input vertices to the output vertices, represents the structure of a computation.
Vertices are decorated by elementary operations on partial recursive functions and 
edges are decorated by partial recursive functions that are inputs and outputs of the vertex
operations, see \cite{Man1}, \cite{Man2}, and see also the discussion in \cite{DelMar}
on generalizations of Manin's Hopf algebra of flow charts.
The computation associated to a graph $\Gamma$ depends on a set of parameters. 

\smallskip

We consider min-plus characters $\psi: \cH \to \bS$, where the choice of the target
min-plus semiring $\bS$ accounts for the dependence on parameters. Typical such
characters would be the running time of the computation (if computations associated
to different connected components of the graph are run sequentially) or the memory size
involved in the computation, with $\psi(\Gamma)=\psi(\Gamma_1)+\psi(\Gamma_2)=
\psi(\Gamma_1)\odot \psi(\Gamma_2)$, for a disjoint union 
$\Gamma = \Gamma_1  \sqcup \Gamma_2$.

\smallskip

In the theory of computation, one approach to characterize the complexity of computable
functions in a machine-independent way is by considering a sequence of machines in
a given class ($1$-tape machines, multiple tape machines, etc.) and associate to each
machine in the sequence a step-counting function, which is the number of steps of tape
(or computing time) that the machine takes to compute a given recursive function
(or infinity if the computation does not stop). This method is the basis for speed-up and
compression theorems, see \cite{Blum} for more details.

\smallskip

Suppose given a decorated graph $\Gamma \in \cH$ with decorations by recursive functions and
operations as in Manin's Hopf algebra of flow charts. Then, given a class of machines,
we let $\psi_n(\Gamma)$ be the step-counting function of the $n$-th machine in
the class, when it computes the output of $\Gamma$. We set $\psi_n(\Gamma)=\infty$
if the $n$-th machine does not halt when fed the input of $\Gamma$. We also assume
that, if $\Gamma$ has several components, the computations are done sequentially,
so that $\psi_n(\Gamma_1 \sqcup \Gamma_2)=\psi_n(\Gamma_1)+\psi_n(\Gamma_2)$,
hence all the $\psi_n$ are $\bT$-valued min-plus characters. Moreover, in the case
of a union that is not disjoint, one can assume that the step-counting functions $\psi_n$
satisfy an inclusion-exclusion principle $\psi_n(\Gamma_1 \cup \Gamma_2)=
\psi_n(\Gamma_1)+\psi_n(\Gamma_2)- \psi_n(\Gamma_1\cap \Gamma_2)$.

\smallskip

We then consider the Rota--Baxter operator of weight $+1$ given by the partial
sum, as in \S \ref{parsumSec}. The preparation of the character 
$\psi(\Gamma)=(\psi_n(\Gamma))_{n\in \N}$ is given by
$$ \tilde\psi_n(\Gamma) = \min \{ \psi_n(\Gamma), \psi_n(\Gamma/\gamma) +
\sum_{k=1}^{n-1} \tilde\psi_k(\gamma) \}, $$
where the minimum is taken over subgraphs $\gamma$.
Notice that, by the inclusion-exclusion property, we are comparing the
size (number of steps/computing time) of $\psi_n(\Gamma)=\psi_n(\Gamma/\gamma)+
\psi_n(\gamma)-\psi_n(\partial\gamma)$, with the size of 
$\psi_n(\Gamma/\gamma)+\sum_{k=1}^{n-1} \psi_k(\gamma)$; and then the minimum
of these with the further terms $\sum_{k=1}^{n-1} (\psi_k(\gamma')+\psi_k(\gamma/\gamma'))$,
for subgraphs $\gamma'\subset \gamma$, and so on, in the recursive structure of the 
$\tilde\psi_k(\gamma)$. At each step, one identifies smaller graphs inside $\Gamma$ for which 
either the cumulative computational time of all the previous machines in
the series is small, or the additional computational cost of the ``interior" part of the subgraph
$\gamma\smallsetminus \partial\gamma$ is small. 

\smallskip

In the case of a graph $\Gamma$ for which the $n$-th machine does not halt, so 
$\psi_n(\Gamma)=\infty$,  the character $\tilde\psi_n$ can be finite, provided 
the following conditions are realized:
\begin{itemize}
\item The source of the infinite computational time for the $n$-th machine was 
localized in an area $\gamma\smallsetminus \partial\gamma$ of the graph
$\Gamma$, that is, $\psi_n(\Gamma/\gamma)<\infty$.
\item None of the previous machines had infinite computational time on this
region of the graph: $\psi_k(\gamma)<\infty$ for all $k=1,\ldots, n-1$.
\end{itemize}

\medskip
\subsection{Nearest neighbor potentials and Markov random fields}

Given a subgraph $\gamma \subset \Gamma$ we denote by $\partial\gamma$ the
subgraph with $E(\partial\gamma)$ the set of edges in $E(\Gamma)$ with
$\partial e$ consisting of a vertex in $V(\gamma)$ and a vertex in
$V(\Gamma)\smallsetminus V(\gamma)$. The set of vertices $V(\partial\gamma)$
is the union of these endpoints, for all $e\in E(\partial\gamma)$.
In particular, for a vertex $v\in V(\Gamma)$ we write $\partial(v)$ for the set 
of vertices $V(\partial\{ v \}) \subset V(\Gamma)$.

\smallskip

A {\em Markov random field} on a graph is a map 
$\pi: \cP(V(\Gamma)) \to \R$, where $\cP(V(\Gamma))$ is the set of subsets of $V(\Gamma)$,
satisfying $\pi(A) >0$ for all $A\in \cP(V(\Gamma))$ and 
\begin{equation}\label{Markovfield}
\frac{\pi(A\cup \{ v \})}{\pi(A)} = \frac{\pi(A \cap \partial(v)) \cup \{ v \})}{\pi(A \cap \partial(v))},
\end{equation}
for all $A\in \cP(V(\Gamma))$ and all $v\in V(\Gamma)$, see \S 1 of \cite{Pres}.

\smallskip

A {\em nearest neighbor potential} on a graph is a function $\cW: \cP(V(\Gamma)) \to \R$
satisfying
\begin{equation}\label{Vnearneigh}
\cW(A \cup \{ v \})- \cW(A) = \cW(A \cap \partial(v)) \cup \{ v \}) - \cW(A \cap \partial(v)),
\end{equation}
for all $A\in \cP(V(\Gamma))$ and all $v\in V(\Gamma)$, see \S 1 of \cite{Pres}.
Unlike \cite{Pres}, here we do not require normalizations for $\pi$ by $\pi(\emptyset)$,
or of $\cW$, by the partition function $Z=\sum_A \exp(\cW(A))$.

\smallskip

We extend the notion of nearest neighbor potentials and Markov random fields from
a single graph to a (finite or infinite) family of graphs. 

\begin{defn}\label{Gmarkov}
Given a family $\cG$  of finite graphs a Markov random field on $\cG$ is a function
$\pi: \cG \to \R$ satisfying $\pi(\Gamma)>0$ for all $\Gamma \subset \cG$ and
\begin{equation}\label{Markovfield2}
\frac{\pi(\Gamma \cup \{ v \})}{\pi(\Gamma)} = \frac{\pi(\Gamma \cap \partial(v)) \cup \{ v \})}{\pi(\Gamma \cap \partial(v))},
\end{equation}
whenever the graphs $\Gamma \cup \{ v \}$, $\Gamma \cap \partial(v)$ and $\Gamma \cap \partial(v)) \cup \{ v \}$
belong to $\cG$. A nearest neighbor potential on $\cG$ is a function $\cW: \cG \to \R$ satisfying
\begin{equation}\label{Vnearneigh2}
\cW(\Gamma \cup \{ v \})- \cW(\Gamma) = \cW(\Gamma \cap \partial(v)) \cup \{ v \}) - \cW(\Gamma \cap \partial(v)),
\end{equation}
whenever $\Gamma \cup \{ v \}$, $\Gamma \cap \partial(v)$ and $\Gamma \cap \partial(v)) \cup \{ v \}$
belong to $\cG$.
\end{defn}

\smallskip

We recover the usual notion of \cite{Pres} if we fix a graph $\Gamma$ and we define
$\cG$ to be the set of all induced subgraphs of $\Gamma$, namely all subgraphs determined
by a choice of a subset $A$ of vertices of $\Gamma$, and all the edges of $\Gamma$ between 
those vertices.

\smallskip

\begin{lem}\label{betaMarkov}
If $\cW: \cG \to \R$ is a nearest neighbor potential on $\cG$, then, for all $\beta>0$, setting 
$\pi_\beta(\Gamma)=e^{-\beta \cW(\Gamma)}$ defines a random Markov field $\pi_\beta:\cG \to \R$.
\end{lem}

\proof This follows immediately by adapting the general observation of 
\S 1 of \cite{Pres},  that if $\pi$ is a Markov random field then $\cW(A)=\log(\pi(A))$
is a nearest neighbor potential and, conversely, given a nearest neighbor potential $\cW$,
setting $\pi(A)=\exp(\cW(A))$ gives a Markov random field. 
\endproof

\smallskip

\smallskip

Let $\cG$ be a family of finite graphs, closed under disjoint unions, 
and $\cH=\cH(\cG)$ the Hopf algebra generated as a commutative algebra
by the connected components of elements of $\cG$ with coproduct 
$\Delta(\Gamma)=\Gamma\otimes 1 + 1\otimes\Gamma + \sum \gamma \otimes \Gamma/\gamma$
where the non-primitive part of the coproduct is the sum over all pairs of a subgraph $\gamma$
and the quotient graph $\Gamma/\gamma$ (where every component of $\gamma$ is contracted to a vertex)
such that both $\gamma$ and $\Gamma/\gamma$ belong to $\cG$.

\smallskip

\begin{lem}\label{tauMarkov}
For $\cG$ and $\cH=\cH(\cG)$ as above, a nearest-neighbor potential $\cW:\cG \to \R$
defines a min-plus character $\cW: \cH \to \bT$.
\end{lem}

\proof It is immediate to check that \eqref{Vnearneigh2} implies
$\cW(\Gamma)=\cW(\Gamma_1) + \cW(\Gamma_2)$ for a disjoint union 
$\Gamma=\Gamma_1 \sqcup \Gamma_2$. In fact, by inductively adding
vertices of $\Gamma_2$ we have $\cW(\Gamma_1 \cup \{ v\}) - \cW(\Gamma_1)=\cW(\{ v \})$
for $v\in V(\Gamma_2)$ and assuming that for $\gamma\subset \Gamma_2$ with
$\# V(\gamma)=n$ we have $\cW(\Gamma_1 \cup \gamma) = \cW(\Gamma_1) + \cW(\gamma)$
we obtain $\cW(\Gamma_1 \cup \gamma\cup \{ v \})=\cW(\Gamma_1) + \cW(\gamma)
+\cW((\gamma \cap \partial(v))\cup \{ v\}) -\cW(\gamma \cap \partial(v)) =
\cW(\Gamma_1) + \cW(\gamma) + \cW(\gamma \cup \{ v \})-\cW(\gamma)=\cW(\Gamma_1) + \cW(\gamma \cup \{ v \})$,
for all $v\in V(\Gamma_2)\smallsetminus V(\gamma)$.
\endproof

\smallskip

Similarly, for more general min-plus semirings $\bS$, 
min-plus characters $\psi: \cH \to \bS$ that depend only on the vertex set of graphs
define $\bS$-valued nearest neighbor potentials.

\begin{lem}\label{Snearneigh}
Let $\bS$ be a min-plus semiring, and 
let $\psi: \cH \to \bS$ be a min-plus character with the property that the value $\psi(\Gamma)\in \bS$
depends only on the set $V(\Gamma)$ of vertices of $\Gamma$. Then $\psi$ is a $\bS$-valued
nearest neighbor potential.
\end{lem}

\proof We have $\psi(\Gamma)=\psi((\Gamma\cap \partial(v)) \cup (\Gamma\smallsetminus \partial(v)))$.
If the value of $\psi$ only depends on the vertex set, then the latter is equal to
$\psi((\Gamma\cap \partial(v)) \sqcup (\Gamma\smallsetminus \partial(v)))$. Since $\psi$ is a min-plus
character, this is $\psi(\Gamma\cap \partial(v)) +\psi (\Gamma\smallsetminus \partial(v))$. Thus, we have
$\psi(\Gamma \cup \{ v \})-\psi(\Gamma)= \psi((\Gamma \cap\partial{v})\cup \{ v\})+ \psi(\Gamma \smallsetminus \partial(v))
- \psi(\Gamma\cap \partial(v)) - \psi (\Gamma\smallsetminus \partial(v)) = \psi((\Gamma \cap\partial{v})\cup \{ v\})
- \psi (\Gamma\smallsetminus \partial(v))$.
\endproof

\smallskip

We can then view the Birkhoff factorization of min-plus characters in thermodynamic semirings as
a method for generating new Markov random fields from given ones.

\begin{prop}\label{MrndfieldBirk}
Let $\cW: \cH \to \bS$ be a nearest neighbor potential, with associated Markov random field 
$\pi_\beta(\Gamma)=e^{-\beta \cW(\Gamma)}$. 
Let $\bS_{\beta,S}$ be the thermodynamic deformation of $\bS$ with $S$ the
Shannon entropy, and let 
Let $T: \bS_{\beta,S} \to \bS_{\beta,S}$ be an $\otimes_{\beta,S}$-linear 
weight-one Rota--Baxter operator. Let $\cW_{\beta,\pm}: \cH\to \bS_{\beta,S}$
be the two parts of the Birkhoff factorization of $\cW$. Then  
$\pi_{\beta,\pm}(\Gamma)=e^{-\beta \cW_\pm(\Gamma)}$ are Markov random fields.
\end{prop}

\proof The factorization is given by
$$ \cW_{\beta,-}(\Gamma) = -\beta^{-1} \log\left( e^{-\beta T \cW(\Gamma)} + \sum e^{-\beta T(\cW_{\beta,-}(\gamma) + \cW(\Gamma/\gamma)) } \right), $$
$$ \cW_{\beta,+}(\Gamma)=-\beta^{-1} \log\left( e^{-\beta \cW_{\beta,-}(\Gamma)} + e^{-\beta \tilde\cW_\beta(\Gamma)}\right),$$
where 
$$ \tilde\cW_\beta(\Gamma)=-\beta^{-1} \log\left( e^{-\beta\cW(\Gamma)} + \sum e^{-\beta (\cW_{\beta,-}(\gamma) + \cW(\Gamma/\gamma)) } \right). $$
According to Theorem \ref{Birk1betaProp}, $\cW_{\beta,\pm}: \cH \to \bS_{\beta,S}$ are min-plus characters.
Moreover, the explicit expression above shows that, if $\cW(\Gamma)$ depends only on the set 
$V(\Gamma)$ of vertices of $\Gamma$, then so do also the $\cW_{\beta,\pm}(\Gamma)$. Thus, by Lemma \ref{Snearneigh}
the $\cW_{\beta,\pm}$ are $\bS_{\beta,S}$-valued nearest neighbor potentials, and $\pi_{\beta,\pm}$
are Markov random fields.
\endproof

\medskip

\subsection{Algebro-geometric Feynman rules and polynomial countability}

In perturbative quantum field theory, one can write the Feynman integrals in the parametric
form as (unrenormalized) period integrals on the complement of the (affine) graph hypersurface
$X_\Gamma \subset \A^{\# E(\Gamma)}$, defined by the vanishing of the graph polynomial
$\Psi_\Gamma(t)=\sum_T \prod_{e\notin E(T)} t_e$, where the sum is over spanning trees
of the Feynman graph $\Gamma$ and $t=(t_e)_{e\in E(\Gamma)}\in \A^{\# E(\Gamma)}$,
see \cite{Mar} for a general overview.

\smallskip

It was observed in \cite{AluMa} that the class in the Grothendieck ring of the affine hypersurface
complement $Y_\Gamma := \A^{\# E(\Gamma)}\smallsetminus X_\Gamma$ determines a 
morphism of commutative rings, from the Hopf algebra of Feynman graphs to the Grothendieck ring,
since it satisfies
\begin{equation}\label{YGamma}
 [Y_\Gamma] = [Y_{\Gamma_1}] \cdot [Y_{\Gamma_2}] 
\end{equation} 
when $\Gamma$ is a disjoint union $\Gamma = \Gamma_1 \sqcup \Gamma_2$. Such
morphisms were termed ``algebro-geometric Feynman rules" in \cite{AluMa}, where examples
based on Chern classes of singular varieties were also constructed, with values in a suitable
Grothendieck group of immersed conical varieties.

\smallskip

Recall that a variety $X$ defines over $\Z$ is polynomially countable if for all the
mod $p$ reductions $X_p$, the counting functions of points over $\F_q$, with $q=p^r$, 
is a polynomial in $q$ with $\Z$-coefficients, namely $N(X,q):=\# X_p(\F_q)=P_X(q)$. Polynomial
countability is a consequence (and, modulo certain conjectures on motives, equivalent) to the
class in the Grothendieck ring $[X]=P_X(\bL)$ being in the polynomial subring $\Z[\bL]$ 
generated by the Lefschetz motive, and to the motive $\m(X)$ being a mixed Tate motive
over $\Z$. An important question in the ongoing investigations of the relations between
quantum field theory and motives is understanding when (for which Feynman graphs) the varieties $X_\Gamma$ 
(or equivalently $Y_\Gamma$) are mixed Tate motives. This question has attracted a
lot of attention in recent years.

\smallskip

We can define a max-plus character related to the behavior of the counting functions
$\# X_p(\F_q)$ for the graph hypersurface complement, that expresses the question
of their polynomial countability.

\begin{lem}\label{minplusq}
Let $\cH$ be the Hopf algebra of Feynman graphs and let $\psi: \cH \to \bT_{max}$ be
defined by $N(Y_\Gamma ,q)\sim q^{\psi(\Gamma)}$, up to lower order terms in $q$, if $Y_\Gamma$
is polynomially countable and $\psi(\Gamma)= -\infty$ if it is not. Then $\psi$ is a max-plus
character, namely $\psi(xy)=\psi(x)+\psi(y)$.
\end{lem}

\proof The counting function $N(X,q)=N([X],q)$ factors through the Grothendieck ring, hence
by \eqref{YGamma} we have $N(Y_\Gamma,q)=N(Y_{\Gamma_1},q) N(Y_{\Gamma_2},q)$
for a disjoint union $\Gamma = \Gamma_1 \sqcup \Gamma_2$. 
If both are polynomially countable, then the exponents of the leading terms satisfy
$\psi(\Gamma)=\psi(\Gamma_1)+\psi(\Gamma_2)$. If at least one of them is not
polynomially countable then $\psi(\Gamma)=-\infty$, which is also equal to $\psi(\Gamma_1)+\psi(\Gamma_2)$,
since one of these terms is also $-\infty$. Thus, the result follows.
\endproof

\medskip

The simplest possible Rota--Baxter operator of weight $-1$ is the identity, $T=id$, which obviously satisfies 
the linearity hypothesis $T(a+b)=T(a)+T(b)$ discussed in \S \ref{minoneSec}. 
Observe that, in the case where the operator $T$ is linear (in the ordinary sense), 
the argument of Proposition \ref{algrensemimin1}
goes through unchanged, if we replace the min-plus tropical semiring
with the analogous max-plus $\bT_{max}$ semiring, by simply replacing
$\oplus=\min$ with $\oplus=\max$, and $\infty$ with $-\infty$ as the additive unit.

\smallskip

In the case of the max-plus character of Lemma \ref{minplusq}, the 
preparation $\tilde\psi(x)$, with respect to $T=id$ then acquires a very simple geometric 
meaning. We have
$$ \tilde\psi(\Gamma)=\max \{ \psi(\Gamma), \tilde\psi(\gamma) + \psi(\Gamma/\gamma) \} = \max\{ \psi(\Gamma),
\sum_{j=1}^N \psi(\gamma_j) + \psi(\gamma_{j-1}/\gamma_j) \} $$
where the maximum is taken over all nested families of subgraphs $\gamma_N \subset \gamma_{N-1}\subset\cdots
\subset \gamma_0=\Gamma$. In the case of a graph $\Gamma$ for which $Y_\Gamma$ is not polynomially countable,
the preparation $\tilde\psi(\Gamma)$ extracts chains of subgraphs and quotient graphs that are polynomially
countable.

\smallskip

A similar example is obtained by considering cases where the graph hypersurfaces $X_\Gamma$ and $Y_\Gamma$
depend on parameters (for example, if one works in the massive, instead of massless case, or if one considers the
hypersurface defined by the second Symanzik polynomial, instead of the first, so that one has the dependence on
the external momenta. In such cases the max-plus character $\psi(\Gamma)$ defined above takes values
in a semiring $\bS$ of functions on the set of parameters, with values in $\bT_{max}$. One can then consider
Rota--Baxter operators of weight $-1$ given by multiplication by the characteristic function of certain subsets
of parameters. The corresponding preparation, as in the simpler case above, would identify 
subgraphs and quotient graphs that are polynomially countable for specific choices of the parameters.

\bigskip
\bigskip
\bigskip

\noindent {\bf Acknowledgment} The first author is supported by NSF grants DMS-1007207, 
DMS-1201512, PHY-1205440. The second author was supported by a Summer Undergraduate 
Research Fellowship at Caltech and by the Rose Hills Foundation.


\begin{thebibliography}{99}

\bibitem{AluMa} P.~Aluffi, M.~Marcolli, {\em Algebro-geometric Feynman rules},
International Journal of Geometric Methods in Modern Physics, Vol.8 (2011) N.1, 203--237.

\bibitem{ArsLeu} M.~Arsove, H.~Leutwiler, {\em Algebraic potential theory},
American Mathematical Society, 1980.

\bibitem{BaFriLei} J.~Baez, T.~Fritz, T.~Leinster, {\em A characterization of entropy 
in terms of information loss}, Entropy 13 (2011), no. 11, 1945--1957. 

\bibitem{BenZyc} I.~Bengtsson, K.~Zyczkowski, {\em Geometry of quantum states}, Cambridge University Press, 2006.

\bibitem{Blum} M.~Blum, {\em A machine-independent theory of the complexity of recursive functions},
Journal of the Association for Computing Machinery, Vol.14 (1967) N.2, 322--336.

\bibitem{CoCo}  A.~Connes, C.~Consani, {\em From monoids to hyperstructures: in search 
of an absolute arithmetic}, in ``Casimir force, Casimir operators and the Riemann hypothesis", 
147--198, Walter de Gruyter, Berlin, 2010. 

\bibitem{CoKr} A.~Connes, D.~Kreimer, {\em Renormalization in quantum field 
theory and the Riemann-Hilbert problem. I. The Hopf
algebra structure of graphs and the main theorem}, Communications in Mathematical Physics 210 (2000), 
no. 1, 249--273.

\bibitem{CoMa} A.~Connes, M.~Marcolli, {\em Noncommutative geometry,
quantum fields and motives}, Colloquium Publications, Vol.55, American
Mathematical Society, 2008.

\bibitem{DelMar} C.~Delaney, M.~Marcolli, {\em Dyson--Schwinger equations in the theory
of computation}, arXiv:1302.5040, to appear in the volume ``Periods and Motives" (Editors: Luis çlvarez-C\'onsul, 
Jos\'e Ignacio Burgos Gil, Kurusch Ebrahimi-Fard, David A. Ellwood), Clay Institute and AMS.

\bibitem{EFGuo}  K.~Ebrahimi-Fard, L.~Guo, {\em Rota-Baxter algebras in renormalization of perturbative quantum field theory}, in ``Universality and renormalization", 47--105, 
Fields Inst. Commun., 50, American Mathematical Society, 2007.

\bibitem{EFGK} K.~Ebrahimi-Fard, L.~Guo, D.~Kreimer, {\em Integrable renormalization. II. 
The general case}, Ann. Henri Poincar\'e 6 (2005), no. 2, 369--395.

\bibitem{GeTsa} M.~Gell-Mann, C.~Tsallis, {\em Nonextensive Entropy}, Oxford University Press, 2004.

\bibitem{Guo} L.~Guo, {\em An Introduction to Rota-Baxter Algebra}, International Press 
and High Education Press, 2012.

\bibitem{KerSchn} S.~Kerkhoff, F.M.~Schneider, {\em A tropical version of the Gelfand representation},
preprint, MATH-AL-14-2012, Technische Universit\"at Dresden, 2012. 

\bibitem{Man1} Yu.I.~Manin, {\em Renormalization and computation, I: motivation and background},
in ``OPERADS 2009", 181--222, S\'emin. Congr., 26, Soc. Math. France, Paris, 2013.

\bibitem{Man2} Yu.I.~Manin, {\em Renormalization and computation II: time cut-off and the halting problem}, 
 Math. Structures Comput. Sci. 22 (2012), no. 5, 729--751.

\bibitem{Man3} Yu.I.~Manin, {\em Infinities in quantum field theory and in classical computing: renormalization program}, in ``Programs, proofs, processes", 307--316, 
Lecture Notes in Comput. Sci., 6158, Springer, 2010. 

\bibitem{Mar} M.~Marcolli, {\em Feynman motives}, World Scientific, 2010. 

\bibitem{MaThor} M.~Marcolli, R.~Thorngren, {\em Thermodynamic semirings}, 
Journal of Noncommutative Geometry, Vol.8 (2014) N.2,
337--392.

 \bibitem{NeStor} S.~Neshveyev, E.~St\/{o}rmer, {\em Dynamical entropy in operator algebras}, 
 Springer, 2006.
 
 \bibitem{PeCli} Y.~Pesin, V.~Climenhaga, {\em Lectures on fractal geometry and dynamical systems},
 American Mathematical Society, 2009.
 
 \bibitem{Pres} C.J.~Preston, {\em Gibbs states on countable sets}, Cambridge University Press, 1974.
 
 \bibitem{Rama} N.~Ramachandran, {\em Zeta functions, Grothendieck groups, and the Witt ring},
 arXiv:1407.1813.
 
 \bibitem{Viro} O.~Viro, {\em Dequantization of real algebraic geometry on logarithmic paper},
European Congress of Mathematics, Vol. I (Barcelona, 2000), 135--146, 
Progr. Math., 201, Birkh\"auser, 2001.
 
 \bibitem{Wilde} M.~Wilde, {\em Quantum information theory}, Cambridge University Press, 
Cambridge, 2013.

\end{thebibliography}
\end{document}